\documentclass[12pt,a4paper]{article}
\usepackage{latexsym}
\usepackage{amsfonts}
\usepackage{amssymb}
\usepackage{amsmath}
\usepackage{wasysym}
\usepackage[mathscr]{eucal}
\usepackage{theorem}
\usepackage{enumerate}
\usepackage{cite}
\usepackage{url}
\usepackage[dvips]{color}
\usepackage{graphicx}
\usepackage{xypic}
\usepackage[dvips]{color}
\usepackage[all]{xy}

\theoremstyle{plain}
\theorembodyfont{\itshape}
\newtheorem{thm}{Theorem}
\newtheorem{lem}{Lemma}
\newtheorem{prp}{Proposition}

\newtheorem{dfn}{Definition}

\theorembodyfont{\upshape}
\newtheorem{exe}{Example}
\newtheorem{rmk}{Remark}

\newcommand{\proof}{\noindent {\bf Proof:} \hspace{0.1in}}
\newcommand{\qed}{\hfill\mbox{\raggedright $\Box$}\medskip}

\newcommand{\smin}{\,\raisebox{0.06em}{${\scriptstyle \in}$}\,}

\newcommand{\smwedge}{{\scriptstyle \wedge\,}}

\newcommand{\smcirc}{{\scriptstyle \,\circ\,}}

\newcommand{\bwedge}[1]{\raisebox{0.2ex}{${\textstyle \bigwedge}$}
 \ensuremath{^{\raisebox{-0.2ex}{${\scriptstyle #1}$}}\,}}

\newcommand{\dbwedge}[2]{\raisebox{0.2ex}{${\textstyle \bigwedge}$}
 \ensuremath{_{\scriptstyle #1}^{\raisebox{-0.2ex}{${\scriptstyle #2}$}}\,}}

\newcommand{\oast}{\ensuremath{\:\!%
 {\raisebox{0.1ex}{$\scriptscriptstyle \bigcirc$} \hspace{-0.5em}
 \ast}}}
\newcommand{\ostar}{\ensuremath{\:\!%
 {\raisebox{0.05ex}{$\scriptscriptstyle \bigcirc$} \hspace{-0.5em}
 \star}}}

\begin{document}

\title{Lagrangian Distributions and Connections \\
       in Symplectic Geometry}
\author{Michael Forger$\,^1\,$%
        \thanks{Work partially supported by CNPq (Conselho Nacional
                de Desenvolvimento Cient\'{\i}fico e Tecno\-l\'ogico),
                Brazil.}
        ~~and~
        Sandra Z. Yepes$\,^2$
}
\date{\normalsize
      Departamento de Matem\'atica Aplicada, \\
      Instituto de Matem\'atica e Estat\'{\i}stica, \\
      Universidade de S\~ao Paulo, \\
      Caixa Postal 66281, \\
      BR--05314-970~ S\~ao Paulo, S.P., Brazil
}
\footnotetext[1]{\emph{E-mail address:} \textsf{forger@ime.usp.br}}
\footnotetext[2]{\emph{E-mail address:} \textsf{samaza@ime.usp.br}}
\maketitle

\thispagestyle{empty}

\begin{abstract}
\noindent
 We discuss the interplay between lagrangian distributions and connections
 in symplectic geometry, beginning with the traditional case of symplectic
 manifolds and then passing to the more general context of poly- and
 multisymplectic structures on fiber bundles, which is relevant for
 the  covariant hamiltonian formulation of classical field theory.
 In particular, we generalize Weinstein's tubular neighborhood theorem
 for symplectic manifolds carrying a (simple) lagrangian foliation to
 this situation. In~all cases, the Bott connection, or an appropriately
 extended version thereof, plays a central role.
\end{abstract}
\begin{center}
 {\small AMS Subject Classification (2000):
         53D12 (Primary), 37J05, 70G45 (Secondary)}
\end{center}
\begin{flushright}
 \parbox{12em}{
  \begin{center}
   Universidade de S\~ao Paulo \\
   RT-MAP-1202 \\
   February 2012
  \end{center}
 }
\end{flushright}

\cleardoublepage

\setcounter{page}{1}

\section{Introduction}

In symplectic geometry, lagrangian foliations or, more generally,
lagrangian distributions (i.e., lagrangian subbundles of the tangent
bundle, which may or may not be involutive) are an important tool.
Although unfamiliar from riemannian or even lorentzian geometry, where
the notion of a lagrangian subbundle does not arise (except in the
rather uninteresting case of two-dimensional lorentzian manifolds),
they also appear in the theory of certain pseudo-riemannian manifolds,
namely those of zero signature.
However, the inter\-play between lagrangian distributions and connections
in symplectic geometry is quite different from, and considerably more
interesting than, in pseudo-riemannian geometry.
In what follows, we shall show that this interplay admits a completely
natural extension to the context of poly- and multisymplectic fiber
bundles, whose precise mathematical definition can be found in~%
\cite{FG} and which appear naturally in the covariant hamiltonian
formulation of classical field theory.

The paper is divided into two parts.
The first part (Sect.~2-5) discusses connections which are compatible
with a given foliation, whereas the second part (Sect.~6-9) uses them
to derive various structure theorems.
More specifically, we begin by reviewing some standard issues from
symplectic geometry: the construction of the Bott connection in~%
Sect.~2 and the classification of symplectic connections in Sect.~3\,:
here, we also prove an analogous classification theorem for symplectic
connections preserving a given lagrangian foliation which~-- although
of independent interest~-- does not seem to have been explicitly
formulated in the literature.
Next, we show how to extend this classification to polysymplectic
fiber bundles in Sect.~4 and to multisymplectic fiber bundles in Sect.~5.
The second part starts, in Sect.~6, with an exposition of
the program to be developed in the remainder of the paper,
followed by a study, in Sect.~7, of manifolds equipped with a
given foliation by flat affine manifolds (where ``flat affine''
refers to a given partial connection along the leaves), since
this is the situation prevailing in all cases of interest here.
The results are applied in Sect.~8 to symplectic manifolds with a
lagrangian foliation, allowing us to give a simple proof, as well as
a generalization, of Weinstein's tubular neighborhood theorem~\cite{Th,We}.
Exactly the same technique leads to what we call the structure theorem
for polysymplectic and multisymplectic fiber bundles, presented in
Sect.~9: it can be viewed as an analogue of the theorem from
symplectic geometry that characterizes which symplectic manifolds
are cotangent bundles (or ``pieces'' of cotangent bundles, possibly
up to coverings).
Finally, in Sect.~10, we present our conclusions.

\section{Lagrangian distributions and the Bott connection}

In this section, we briefly review the definition of the Bott connection
for symplectic mani\-folds carrying a given lagrangian foliation, noting
that a completely analogous concept also exists for pseudo-riemannian
manifolds of zero signature (these are the only ones whose tangent
spaces admit lagrangian subspaces).

\pagebreak

As a preliminary step, we recall that given a manifold~$M$, a vector bundle
$V$ over~$M$ and a distribution~$L$ on~$M$, a \emph{partial linear connection
in $V$ along~$L$} is an $\mathbb{R}$-bilinear map\,%
\footnote{Given a manifold~$M$ and a vector bundle $W$ over~$M$, we denote by
$\mathfrak{F}(M)$ the algebra of (smooth) functions on~$M$ and by $\Gamma(W)$
the space (and $\mathfrak{F}(M)$-module) of (smooth) sections of~$W$. If $W$
is a distribution on~$M$, i.e., a vector subbundle of the tangent bundle $TM$
of~$M$, we shall use the intuitively more appealing expression ``vector field
on~$M$ along~$W$'' for a section of~$W$.}
\begin{equation} \label{eq:PARCON}
 \begin{array}{cccc}
  \nabla: & \Gamma(L) \times \Gamma(V) & \longrightarrow &     \Gamma(V)
  \\[1mm]
          &           (X,s)            &   \longmapsto   & \nabla_{\!X}^{} s
 \end{array}
\end{equation}
which is $\mathfrak{F}(M)$-linear in~$X$ and a derivation in~$s$, i.e.,
satisfies the usual Leibniz rule
\begin{equation} \label{eq:LEIBNR}
 \nabla_{\!X}^{} \, (fs)~=~f \, \nabla_{\!X}^{} s \, + \, (X \cdot f) \, s
 \qquad \mbox{for $\, X \in \Gamma(L) \,$, $f \in \mathfrak{F}(M) \,$,
                  $s \in \Gamma(V)$}~.
\end{equation}
Of course, this gives back the usual definition of a ``full'' linear
connection when $\, L = TM$ \linebreak and hence $\Gamma(L)$ is the
Lie algebra $\mathfrak{X}(M)$ of all (smooth) vector fields on~$M$.
Clearly, the usual definitions of curvature and, in the special case
when $\, V=L$, of torsion also work for partial linear connections if
$L$ is supposed to be involutive.

Obviously, partial linear connections can be obtained from ``full''
ones by \emph{restriction}\/, that is, by restricting the definition
of the covariant derivative of a section from general vector fields 
on~$M$ to vector fields on~$M$ along a given vector subbundle $L$
of~$TM$, and conversely, one may ask whether a given partial linear
connection admits an \emph{extension}\/ to a ``full'' one (from which
it can be derived by restriction), and if so, how one can classify
all possible extensions.

A particularly nice example of a partial linear connection which is
not evidently the restriction of a ``full'' one is the \emph{Bott
connection}\/ associated with any involutive distribution on any
manifold~$M$: denoting by $L^\perp$ the \emph{annihilator}\/ of~$L$,
which by definition is the vector subbundle of the cotangent bundle
$T^* M$ of~$M$ consisting of $1$-forms that vanish on~$L$, this is a
partial linear connection $\nabla^B$ in $L^\perp$ along~$L$ defined by
\begin{equation} \label{eq:BOTTC1}
 \nabla_{\!X}^B \alpha~=~\mathbb{L}_X^{} \alpha
 \qquad \mbox{for $\, X \in \Gamma(L)$, $\alpha \in \Gamma(L^\perp)$}~,
\end{equation}
where $\mathbb{L}_X^{}$ denotes the Lie derivative (of $1$-forms) along~$X$,
or more explicitly,
\begin{equation} \label{eq:BOTTC2}
 \begin{array}{c}
  (\nabla_{\!X}^B \alpha)(Y)~=~X \!\cdot \alpha(Y) \, - \, \alpha([X,Y]) \\[2mm]
  \mbox{for $\, X \in \Gamma(L)$, $\alpha \in \Gamma(L^\perp)$,
            $Y \in \mathfrak{X}(M)$}~.
 \end{array}
\end{equation}
To show that this is really a partial linear connection in~$L^\perp$ along~$L$,
suppose that $X$ is a vector field on~$M$ along~$L$ and note that (a)~for any
section $\alpha$ of~$L^\perp$ and any vector field $Y$ on~$M$ along~$L$, the
expression in eqn~(\ref{eq:BOTTC2}) vanishes because $L$ is supposed to be
involutive, which means that the Lie derivative along~$X$ really does map
sections $\alpha$ of~$L^\perp$ to sections $\mathbb{L}_X^{} \alpha$ of~%
$L^\perp$, and (b)~for any function~$f$ on~$M$, any section $\alpha$ of~%
$L^\perp$ and any vector field $Y$ on~$M$ (not necessarily along~$L$), we have
\begin{eqnarray*}
 (\mathbb{L}_{fX}^{} \alpha)(Y) \!\!
 &=&\!\! (fX) \cdot \alpha(Y) - \alpha([fX,Y]) \\
 &=&\!\! f (X \cdot \alpha(Y)) - f \alpha([X,Y]) + (Y \cdot f) \, \alpha(X) \\
 &=&\!\! f (\mathbb{L}_X^{} \alpha)(Y) + (Y \cdot f) \, \alpha(X)~,
\end{eqnarray*}
and the last term vanishes because $\alpha$ annihilates $X$, so
$\nabla_X^B$ is really $\mathfrak{F}(M)$-linear in~$X$. \linebreak
Of course, the Bott connection is \emph{flat}\/: its curvature vanishes
trivially, due to the definition of the Lie bracket of vector fields
by means of the formula $\; [\mathbb{L}_X^{},\mathbb{L}_Y^{}]
= \mathbb{L}_{[X,Y]} \,$.

Now assume that $M$ is a manifold of even dimension $2n$ and
$L$ is an involutive distribution on~$M$ which is lagrangian
% either with respect to a given pseudo-riemannian metric~$g$
% of~zero signature or
with respect to a given almost symplectic form (i.e., non-%
degenerate $2$-form)~$\,\omega$.
Then the ``musical isomorphism''~\cite[p.~166]{AM}
% \begin{equation} \label{eq:MUSISO1}
%  g^\flat: TM~\longrightarrow~T^* M \qquad \mbox{with inverse} \qquad
%  g^\sharp: T^* M~\longrightarrow~TM
% \end{equation}
% or
\begin{equation} \label{eq:MUSISO2}
 \omega^\flat: TM~\longrightarrow~T^* M \qquad \mbox{with inverse} \qquad
 \omega^\sharp: T^* M~\longrightarrow~TM
\end{equation}
restricts to an isomorphism
% \begin{equation} \label{eq:MUSISO3}
%  g^\flat: L~\longrightarrow~L^\perp \qquad \mbox{with inverse} \qquad
%  g^\sharp: L^\perp~\longrightarrow~L
% \end{equation}
% or
\begin{equation} \label{eq:MUSISO4}
 \omega^\flat: L~\longrightarrow~L^\perp \qquad \mbox{with inverse} \qquad
 \omega^\sharp: L^\perp~\longrightarrow~L
\end{equation}
which can be used to transfer the Bott connection as defined previously
to a partial linear connection in~$L$ along~$L$: by abuse of language,
it will simply be called the \emph{Bott connection in $L$}. Explicitly,
it is determined by the formula
% \begin{equation} \label{eq:BCPR}
%  g(\nabla_{\!X}^B Y,Z)~
%  =~X \cdot g(Y,Z) \, - \, g(Y,[X,Z])
%  \qquad \mbox{for $\, X,Y \in \Gamma(L) \,$, $Z \in \mathfrak{X}(M)$}
% \end{equation}
% in the pseudo-riemannian case of zero signature and similarly by the formula
\begin{equation} \label{eq:BCSY1}
 \omega(\nabla_{\!X}^B Y,Z)~
 =~X \cdot \omega(Y,Z) \, - \, \omega(Y,[X,Z])
 \qquad \mbox{for $\, X,Y \in \Gamma(L) \,$, $Z \in \mathfrak{X}(M)$}~.
\end{equation}
% in the almost symplectic case.
This connection has an intuitively appealing interpretation: it is nothing
else than a canonical family of ordinary linear connections in the leaves
of the foliation generated by~$L$.
Moreover, the Bott connection in~$L^\perp$ being flat, so is the Bott
connection in~$L$.
But regarding the latter, we can ask for more: we can ask whether it
also has vanishing torsion, since if so, we may conclude that the
leaves of the foliation generated by~$L$ are flat affine manifolds.
Regarding this question, we have the following simple answer.
\begin{thm}~ \label{thm:TORBCSY}
 Let\/~$M$ be a manifold equipped with an almost symplectic form~$\,\omega$
 and let\/ $L$ be any involutive lagrangian distribution on\/~$M$. Then if
 $\,\omega$ is closed, the Bott connection\/ $\nabla^B$ in\/~$L$ has zero
 torsion. More generally, the torsion tensor\/ $T^B$ of\/~$\nabla^B$ is
 related to the exterior derivative of~$\,\omega$ by the formula
 \begin{equation} \label{eq:TBCSY1}
  d \>\! \omega (X,Y,Z)~=~\omega(T^B(X,Y),Z)
 \qquad \mbox{for $\, X,Y \in \Gamma(L) \,$, $Z \in \mathfrak{X}(M)$}~.
 \end{equation}
\end{thm}
\proof
 Writing out the Cartan formula for the exterior derivative of~$\,\omega$,
 \begin{eqnarray*}
  d \>\! \omega(X,Y,Z) \!
  &=&\! X \!\cdot \omega(Y,Z) \, - \, Y \!\cdot \omega(X,Z) \, + \,
        Z \!\cdot \omega(X,Y) \\[1mm]
  & &\! \mbox{} - \, \omega([X,Y],Z) \, + \, \omega([X,Z],Y) \, - \,
                     \omega([Y,Z],X])~.
 \end{eqnarray*}
 and assuming $X$ and $Y$ to be along~$L$, we see that the third of the
 six terms on the rhs of this equation vanishes since $L$ is supposed to
 be isotropic, so using the definition of the torsion tensor combined
 with eqn~(\ref{eq:BCSY1}) to give
 \begin{eqnarray*}
 \lefteqn{\omega \bigl( T^B(X,Y) \,, Z \bigr)~
          =~\omega \bigl( \nabla_{\!X}^B Y - \nabla_{\!Y}^B X - [X,Y] \,,
                          Z \bigr)} \hspace{5mm} \\[1mm]
  &=&\! X \!\cdot \omega(Y,Z) \, - \, \omega(Y,[X,Z]) \, - \,
        Y \!\cdot \omega(X,Z) \, + \, \omega(X,[Y,Z]) \, - \, \omega([X,Y],Z)~,
 \end{eqnarray*}
 we arrive at eqn~(\ref{eq:TBCSY1}), which proves the remaining statements.
\qed

\pagebreak

\noindent
Of course, this result has been known for a long time; see, e.g., Theorem 7.7
of Ref.~\cite{We}. The only difference is that we propose a more systematical
and ample use of the term ``Bott connection''.%
\footnote{In the case of a lagrangian distribution which is involutive
with respect to a pseudo-riemannian metric~$g$ of zero signature, the
construction is completely analogous, and it can be shown that the Bott
connection $\nabla^B$ coincides with the restriction of the Levi-Civita
connection $\nabla$ if and only if the former has zero torsion; more
generally, the difference between the two is proportional to the
torsion tensor of~$\nabla^B$.}

\section{Symplectic connections}

Given a manifold~$M$ equipped with a symplectic form~$\,\omega$ and an
involutive distribution $L$ on~$M$, we can ask the following question:
is the Bott connection $\nabla^B$ in~$L$ the restriction of some
torsion-free symplectic connection $\nabla$ on~$M$, and if so,
what is the set of such torsion-free symplectic connections?

To gain a better understanding of this question and of its importance for
quantization (geometric quantization as well as deformation quantization),
let us briefly explain a few well-known facts about symplectic connections.
We begin with their definition which~-- even though it is standard~-- will
be stated explicitly in order to clarify the terminology.
\begin{dfn}~ \label{def:SYCON}
 Let\/~$M$ be a manifold equipped with an almost symplectic form~%
 $\,\omega$. A~linear connection\/~$\nabla$ on\/~$M$ is said to be
 a \textbf{symplectic connection} if it preserves~$\,\omega$, i.e.,
 satisfies $\, \nabla \omega = 0$, or explicitly,
 \begin{equation} \label{eq:SYCONN}
  X \!\cdot \omega(Y,Z)~
  =~\omega(\nabla_{\!X}^{} Y,Z) \, + \, \omega(Y,\nabla_{\!X}^{} Z)
  \qquad \mbox{for $\, X,Y,Z \in \mathfrak{X}(M)$}~.
 \end{equation}
 The same terminology is used for partial linear connections.
\end{dfn}
In particular, we do not adhere to the convention adopted by some authors
who incorporate the condition of being torsion-free into the definition
of a symplectic connection.

Regarding the question of whether there exist any torsion-free symplectic
connections at all, we begin by noting the following elementary and well
known proposition, which can be viewed as an analogue of Theorem~%
\ref{thm:TORBCSY} for ``full'' linear connections.
\begin{prp}~ \label{prp:SYMCON1}
 Let\/~$M$ be a manifold equipped with an almost symplectic form~$\,\omega$.
 Then if there exists a torsion-free symplectic connection\/ $\nabla$ on\/~$M$,
 $\omega$ must be closed. More generally, the torsion tensor\/ $T$ of a
 symplectic connection\/~$\nabla$ on\/~$M$ is related to the exterior
 derivative of~$\,\omega$ by the formula
 \begin{equation}
  d \>\! \omega (X,Y,Z)~
  =~\omega(T(X,Y),Z) \, + \, \omega(T(Y,Z),X) \, + \, \omega(T(Z,X),Y)~.
 \end{equation}
\end{prp}
\proof
 This is a special case of Lemma~\ref{lem:CONN2} (eqn~(\ref{eq:LEMC1}))
 in Appendix~A.
\qed

\noindent
Conversely, it is well known that on any symplectic manifold, there exist
torsion-free symplectic connections~\cite{Vey,BFFLS}. An explicit proof
can be found in Sect.~2.1 of Ref.~\cite{BCGRS}: it is based on modifying
a given torsion-free linear connection by adding a judiciously chosen
tensor field in order to arrive at a torsion-free symplectic connection.
However, the method can be easily generalized so as to start out from an
arbitrary linear connection and get what one wants in a single stroke:
\begin{prp}~ \label{prp:SYMCON2}
 Let\/~$M$ be a manifold equipped with a symplectic form~$\,\omega$ and\/
 $\nabla^0$ a general linear connection on\/~$M$ with torsion tensor\/~$T^0$.
 Then the formula
 \begin{equation} \label{eq:SYMCON1}
  \begin{array}{rcl}
   \omega \bigl( \nabla_{\!X}^{} Y , Z \bigr) \!
   &=&\! \omega \bigl( \nabla_{\!X}^0 Y , Z \bigr) \, + \,
         {\textstyle \frac{1}{3}} \, (\nabla_{\!X}^0 \omega)(Y,Z) \, + \,
         {\textstyle \frac{1}{3}} \, (\nabla_{\!Y}^0 \omega)(X,Z) \\[2mm]
   & &\! \mbox{} - \,
         {\textstyle \frac{1}{2}} \; \omega \bigl( T^0(X,Y) , Z \bigr) \, + \,
         {\textstyle \frac{1}{6}} \; \omega \bigl( T^0(Z,X) , Y \bigr) \, + \,
         {\textstyle \frac{1}{6}} \; \omega \bigl( T^0(Z,Y) , X \bigr)
  \end{array}
 \end{equation}
 or equivalently
 \begin{eqnarray} \label{eq:SYMCON2}
   \omega \bigl( \nabla_{\!X}^{} Y , Z \bigr) \!
   &=&\! {\textstyle \frac{1}{6}} \;
         \omega \bigl( \nabla_{\!X}^0 Y , Z \bigr) \, + \,
         {\textstyle \frac{1}{6}} \;
         \omega \bigl( \nabla_{\!Y}^0 X , Z \bigr) \, + \,
         {\textstyle \frac{1}{6}} \;
         \omega \bigl( \nabla_{\!Z}^0 X , Y \bigr) \, + \,
         {\textstyle \frac{1}{6}} \;
         \omega \bigl( \nabla_{\!Z}^0 Y , X \bigr) \nonumber \\[1mm]
   & &\! \mbox{} + \, {\textstyle \frac{1}{6}} \;
         \omega \bigl( \nabla_{\!X}^0 Z , Y \bigr) \, + \,
         {\textstyle \frac{1}{6}} \;
         \omega \bigl( \nabla_{\!Y}^0 Z , X \bigr) \, + \,
         {\textstyle \frac{1}{3}} \, X \!\cdot \omega(Y,Z) \, + \,
         {\textstyle \frac{1}{3}} \, Y \!\cdot \omega(X,Z) \qquad \\[1mm]
   & &\! \mbox{} + \,
         {\textstyle \frac{1}{2}} \; \omega \bigl( [X,Y] , Z \bigr) \, - \,
         {\textstyle \frac{1}{6}} \; \omega \bigl( [Z,X] , Y \bigr) \, - \,
         {\textstyle \frac{1}{6}} \; \omega \bigl( [Z,Y] , X \bigr) \nonumber
 \end{eqnarray}
 defines a torsion-free symplectic connection\/ $\nabla$ on\/~$M$.
\end{prp}

\noindent
Obviously, when $\nabla^0$ is itself torsion-free and symplectic,
then $\, \nabla = \nabla^0$.
In passing, we also note that if $\nabla^0$ is torsion-free but
not symplectic, then eqn~(\ref{eq:SYMCON2}) simplifies to
\begin{equation} \label{eq:SYMCON3}
 \begin{array}{rcl}
  \omega \bigl( \nabla_{\!X}^{} Y , Z \bigr) \!
  &=&\! {\textstyle \frac{2}{3}} \;
        \omega \bigl( \nabla_{\!X}^0 Y , Z \bigr) \, - \,
        {\textstyle \frac{1}{3}} \;
        \omega \bigl( \nabla_{\!Y}^0 X , Z \bigr) \, + \,
        {\textstyle \frac{1}{3}} \;
        \omega \bigl( \nabla_{\!X}^0 Z , Y \bigr) \, + \,
        {\textstyle \frac{1}{6}} \;
        \omega \bigl( \nabla_{\!Y}^0 Z , X \bigr) \\[2mm]
  & &\! \mbox{} + \,
        {\textstyle \frac{1}{3}} \, X \!\cdot \omega(Y,Z) \, + \,
        {\textstyle \frac{1}{3}} \, Y \!\cdot \omega(X,Z)~.
 \end{array}
\end{equation}
whereas if $\nabla^0$ is symplectic but not torsion-free,
then eqn~(\ref{eq:SYMCON2}) simplifies to
\begin{eqnarray} \label{eq:SYMCON4}
 \omega \bigl( \nabla_{\!X}^{} Y , Z \bigr) \!
 &=&\! {\textstyle \frac{1}{2}} \;
       \omega \bigl( \nabla_{\!X}^0 Y , Z \bigr) \, + \,
       {\textstyle \frac{1}{2}} \;
       \omega \bigl( \nabla_{\!Y}^0 X , Z \bigr) \, + \,
       {\textstyle \frac{1}{6}} \;
       \omega \bigl( \nabla_{\!Z}^0 X , Y \bigr) \, + \,
       {\textstyle \frac{1}{6}} \;
       \omega \bigl( \nabla_{\!Z}^0 Y , X \bigr)~~ \quad \nonumber \\[1mm]
 & &\! \mbox{} - \, {\textstyle \frac{1}{6}} \;
       \omega \bigl( \nabla_{\!X}^0 Z , Y \bigr) \, - \,
       {\textstyle \frac{1}{6}} \;
       \omega \bigl( \nabla_{\!Y}^0 Z , X \bigr) \\[1mm]
 & &\! \mbox{} + \,
       {\textstyle \frac{1}{2}} \; \omega \bigl( [X,Y] , Z \bigr) \, - \,
       {\textstyle \frac{1}{6}} \; \omega \bigl( [Z,X] , Y \bigr) \, - \,
       {\textstyle \frac{1}{6}} \; \omega \bigl( [Z,Y] , X \bigr)~. \nonumber
\end{eqnarray}
Note the similarity, but also the differences, between these formulas for
symplectic mani\-folds and the definition of the Levi-Civita connection for
pseudo-riemannian manifolds.
In both cases, existence of torsion-free compatible connections is guaranteed,
but in sharp contrast with the pseudo-riemannian case, torsion-free symplectic
connections are far from unique: rather, one can show that the set of all such
connections constitutes an affine space whose difference vector space can be
identified with the space of all totally symmetric tensor fields of rank~$3$;
see, e.g., Refs~\cite{Vey,BFFLS} and, for an explicit proof, Sect.~2.1 of
Ref.~\cite{BCGRS}.
This ambiguity has important implications in mathematical physics, being
closely related to the famous factor ordering problem of quantum mechanics.
More specifically, there is a famous construction of star products in
deformation quantization~\cite{Fed}, now commonly known as the \emph{Fedosov
construction}, which uses as one of its essential ingredients a torsion-free
symplectic connection on classical phase space, and different choices lead
to different factor ordering rules.
This ground-breaking contribution to the quantization problem has even
led some authors to refer to symplectic manifolds equipped with a fixed
torsion-free symplectic connection as \emph{Fedosov manifolds}~\cite{GRS},
and it provides compelling motivation for geometers to study the question
as to what further restrictions on the choice of torsion-free symplectic
connections are implied by introducing additional covariantly constant
geometric structures~-- ideally to the point of singling out a unique
representative.

Of course, there are many possible such structures, among which
we may mention, as particularly important and interesting examples,
K\"ahler manifolds and hamiltonian $G$-spaces; an overview can be
found in Ref.~\cite{BCGRS}. Here, we shall study a specific type,
given by the choice of a lagrangian distribution. This is the kind
of additional structure one meets in geometric quantization~%
\cite{Woo}, and the question of how to construct torsion-free
sym\-plectic connections compatible with it has first been
investigated by He\ss~\cite{He1,He2}. (Somewhat more generally,
geometric quantization uses lagrangian vector subbundles of the
complexified tangent bundle, called \emph{polarizations}, but
we shall in this paper restrict ourselves to real polarizations,
for the sake of simplicity.) The main theorem of He\ss\ regarding
this question, stated in Ref.~\cite{He1} and proved in detail in 
Ref.~\cite{He2}, states that given two involutive lagrangian
distributions $L_1$ and $L_2$ which are transversal, there is
a \emph{unique} symplectic connection preserving both of them:
it has come to be known as the \emph{bilagrangian} connection.%
\footnote{Actually, the theorem of He\ss\ is significantly more
general because it applies even when the lagrangian distributions
are not involutive or the form~$\omega$ is not closed, but we
shall not go into this here.} \linebreak
Another way of looking at this result
% which however works only when $L_1$ and~$L_2$ are involutive
% and $\,\omega$ is closed,
is in terms of pseudo-riemannian geometry, since in this case the
bilagrangian connection, being torsion-free, is simply the Levi-%
Civita connection associated with the pseudo-riemannian metric $g$
of zero signature that can be constructed naturally from~$\,\omega$
together with $L_1$ and $L_2$~\cite{ES,EST} by setting
\[
 g \bigl( X,Y \bigr)~
 =~\omega \bigl( (\mathrm{pr}_1^{} - \mathrm{pr}_2^{})X,Y \bigr)~,
\]
where $\mathrm{pr}_1^{}$ and $\mathrm{pr}_2^{}$ is the projection
onto $L_1$ along~$L_2$ and onto $L_2$ along~$L_1$, respectively.

However, in order to generalize the construction of adequate symplectic
connections to the poly- and multisymplectic framework, we must focus
on the situation where we are given a single involutive lagrangian
distribution $L$ on~$M$, rather than two transversal ones. \linebreak
The question is whether there always exists a torsion-free symplectic
connection on~$M$ that preserves~$L$ and, if so, what is the affine
space of all such connections.
This is a problem of independent interest even within the traditional
context of symplectic geometry, and one that seems to have received
little attention so far.

As a first step in this direction, we note the following extension of
Proposition~\ref{prp:SYMCON1}.
\begin{prp}~ \label{prp:SYMCON3}
 Let\/~$M$ be a manifold equipped with an almost symplectic form~$\,\omega$
 and let\/ $L$ be a lagrangian distribution on\/~$M$. Then if there exists
 a torsion-free symplectic connection\/ $\nabla$ on\/~$M$ preserving\/~$L$,
 $\omega$ must be closed and\/ $L$ must be involutive. In this case, the
 restriction of any such connection to\/~$L$ coincides with the Bott
 connection in\/~$L$.
\end{prp}

\pagebreak

\begin{rmk}~ \label{rmk:SYMCON}
 The last statement is valid under much less restrictive assumptions on the
 torsion tensor~$T$ of~$\nabla$ than stated above: it suffices that $T(X,Y)$
 should be along~$L$ whenever at least one of its arguments is along~$L$.
\end{rmk}
\proof
 The first statement has been proved in Proposition~\ref{prp:SYMCON1}.
 The second statement follows directly from Lemma~\ref{lem:CONN1} in
 Appendix~A. For the third statement, let us assume that $\nabla$ is any
 symplectic connection on~$M$ preserving~$L$ with torsion tensor~$T$
 such that $T(X,Y)$ is along~$L$ as soon as $X$ or $Y$ is along~$L$.
 Then the claim is equivalent to the condition that for all $\, X,Y
 \in \Gamma(L) \,$ and $\, Z \in \mathfrak{X}(M)$,
 \[
  \omega(\nabla_{\!X}^B Y,Z)~=~\omega(\nabla_{\!X}^{} Y,Z)~,
 \]
 which can be derived by comparing eqn~(\ref{eq:BCSY1}) with eqn~%
 (\ref{eq:SYCONN}) taking into account that
 \[
  \omega(Y,[X,Z])~=~\omega(Y,\nabla_{\!X}^{} Z)
 \]
 since $L$ being isotropic and stable under $\nabla$, the expressions
 $\,\omega(Y,\nabla_{\!Z}^{} X)$ and $\,\omega(Y,T(X,Z))$ vanish under
 these assumptions.
\qed

\noindent
Conversely, we can use a partition of unity argument to prove that the
conditions stated in Proposition~\ref{prp:SYMCON3} ($\omega$ is closed
and $L$ is involutive) are not only necessary but also sufficient to
guarantee existence of torsion-free symplectic connections preserving
the distribution~$L$.
\begin{thm}~ \label{thm:LAGCON}
 Let\/~$M$ be a manifold equipped with a symplectic form~$\,\omega$
 and let\/ $L$ be an involutive lagrangian distribution on\/~$M$.
 Then there exist torsion-free symplectic connections~$\,\nabla$
 on\/~$M$ preserving\/~$L$, and the set of all such connections
 constitutes an affine space whose difference vector space can
 be identified with the space of all tensor fields of rank\/~$3$
 on\/~$M$ which (a) are totally symmetric and (b) vanish when%
 ever at least two of their arguments are along\/~$L$.
\end{thm}
\proof
 Concerning existence, we can under the hypotheses of the theorem apply
 the Darboux theorem to guarantee that locally (i.e., on a sufficiently
 small open neighborhood of each point of~$M$), there exists a torsion-%
 free flat symplectic connection on~$M$ preserving~$L\,$: it is simply the
 linear connection on~$M$ whose Christoffel symbols vanish identically in
 these coordinates. (Here, we use a strengthened version of the Darboux
 \linebreak theorem which guarantees the existence of a system of local
 coordinates $(q^i,p\>\!_i^{})$ around each point such that not only
 $\,\omega$ takes the standard form $\, dq^i \,\smwedge\, dp_i^{} \,$
 but also $L$ is generated by the $\partial/\partial p_i^{}$, say; a
 detailed proof can be found, for example, in~\cite[Theorem~1.1]{VA}.)
 Now using a covering of~$M$ by such Darboux coordinate neighborhoods,
 passing to a \linebreak locally finite refinement $(U_\alpha)_%
 {\alpha \in A}$, denoting the corresponding family of linear
 connections by $(\nabla_{\!\alpha})_{\alpha \in A}$ and choosing
 a partition of unity $(\chi_\alpha)_{\alpha \in A}$ subordinate
 to the open covering $(U_\alpha)_{\alpha \in A}$, we can define
 \[
  \nabla~=~\sum_{\alpha \in A} \, \chi_\alpha \, \nabla_{\!\alpha}~.
 \]
 Then it is clear that $\nabla$ preserves~$\,\omega$ as well as~$L$ and
 is torsion-free, since this is true for each $\nabla_{\!\alpha}$ and since
 the conditions of preserving a given differential form, of preserving a
 given vector subbundle and of being torsion-free are all local (i.e.,
 behave naturally under restriction to open subsets) as well as affine.%
 \footnote{The situation with respect to curvature is different because the
 condition of being flat, although still local, is not affine, so although
 each $\nabla_{\!\alpha}$ is flat, this will in general no longer be true
 for~$\nabla$. However, the curvature of~$\nabla$ does vanish when evaluated
 on two vector fields along~$L$, since there $\nabla$ coincides with the Bott
 connection, which is flat.\label{fn:CURV}}
 Regarding uniqueness, or rather the amount of non-uniqueness, we can write
 the difference between any linear connection $\nabla'$ on~$M$ and a fixed
 torsion-free symplectic connection~$\nabla$ on~$M$ preserving~$L$ in the form
 \[
  \nabla_{\!X}' Y~=~\nabla_{\!X}^{} Y \, + \, S(X,Y)~.
 \]
 Moreover, we introduce a (covariant) tensor field $\,\omega_S^{}$ of rank~$3$
 on~$M$ which, due to non-degeneracy of~$\,\omega$, carries exactly the same
 information as~$S$ itself, given by
 \[
  \omega_S^{}(X,Y,Z)~=~\omega(S(X,Y),Z)~.
 \]
 Then it is clear that $\nabla'$ will be torsion-free if and only if $S$
 is symmetric, or equivalently, $\omega_S^{}$ is symmetric in its first
 two arguments, that $\nabla'$ will be symplectic if and only if $S$
 satisfies the identity
 \[
  \omega(S(X,Y),Z) \, + \, \omega(Y,S(X,Z))~=~0~,
 \]
 or equivalently, $\omega_S^{}$ is symmetric in its last two arguments,
 and that $\nabla'$ will preserve~$L$ if and only if $S(X,Y)$ is along~$L$
 whenever $X$ or $Y$ is along~$L$, or equivalently, $\omega_S^{}$ vanishes
 whenever at least two of its arguments are along~$L$.
\qed

These considerations show that there are (at least) three rather different
methods for proving existence of torsion-free symplectic connections:
(a) by modifying a given linear connection through addition of an
appropriately chosen tensor field (see Proposition~\ref{prp:SYMCON2}),
(b) by employing the construction of the bilagrangian connection due
to He\ss\ and (c) by a partition of unity argument.
For our purposes, however, the first two are not fully adequate
since the first provides connections that may not preserve any lagrangian
distribution (this is not enough), whereas the second provides connections
preserving two transversal lagrangian distributions (this is too much).
The problem with the second construction is that the bilagrangian connection
associated with two lagrangian distributions is torsion-free if and only if
both of them are involutive. Therefore, it must be modified when one wants
to deal with situations where one is given a naturally defined involutive
lagrangian distribution~$L$ which has no distinguished lagrangian complement
and which may not even admit any involutive lagrangian complement at all:
an important example is provided by cotangent bundles where the lagrangian
foliation given by the structure as a vector bundle admits tranversal
lagrangian submanifolds (such as the zero section or, more generally,
the graph of any closed $1$-form) but no natural transversal lagrangian
foliation. Such a modification can always be performed by applying the
construction of Proposition~\ref{prp:SYMCON2} to the bilagrangian connection
associated with an arbitrarily chosen lagrangian complement $L'$ of~$L$,
which preserves $L'$ but has non-vanishing torsion (except when $L'$ is
involutive), trading it for what we might call a lagrangian connection,
which no longer preserves $L'$ (except when $L'$ is involutive) but has
vanishing torsion. However, this procedure is somewhat artificial, and
as it turns out, it cannot be extended to the poly- and multisymplectic
setting to be discussed in the next two sections~-- in contrast to the
method based on a partition of unity argument, which extends in a
completely straightforward manner.

\section{Polysymplectic connections}

We begin by stating the definition of a polysymplectic structure,
as given in~Ref.~\cite{FG}.
To this end, we recall first of all that given a fiber bundle~$P$ over
a manifold $M$ with bundle projection~$\; \pi: P \longrightarrow M \;$
and corresponding vertical bundle~$V\!P$ (the  kernel of the tangent
map~$\; T\pi: TP \longrightarrow TM$ of~$\pi$), a vertical vector
field on~$P$ is a section of~$V\!P$ whereas  a vertical $r$-form
on~$P$ is a section of the $r$-th exterior power of the dual bundle
$V^*P$ of~$V\!P$ and, more generally, a (totally covariant) vertical
tensor field of rank~$r$ on~$P$ is a section of the $r$-th tensor
power of the dual bundle $V^*P$ of~$V\!P$.
Similarly, given an additional auxiliary vector bundle~$\hat{T}$
over~$M$, a vertical $r$-form on~$P$ and, more generally, a (totally
covariant) vertical tensor field of rank~$r$ on~$P$ taking values
in~$\hat{T}$~-- or more precisely, in the pull-back $\pi^*(\hat{T})$
of~$\hat{T}$ to~$P$~-- is a section of the tensor product of the
aforementioned exterior power\,/\,tensor power with~$\pi^*(\hat{T})$.
In what follows, we shall denote the Lie algebra of vertical vector
fields on~$P$ by $\mathfrak{X}_V^{}(P)$ and the space of vertical
$r$-forms on~$P$ taking values in~$\hat{T}$ by~$\Omega_V^{\,r}
(P;\pi^* \hat{T})$.% 
\footnote{Note that speaking of vertical forms or (totally covariant)
tensor fields constitutes a certain abuse of language because these
are really equivalence classes of ordinary differential forms or
(totally covariant) tensor fields.} \linebreak
For such forms, there is a complete Cartan calculus, strictly analogous to
the Cartan calculus for (vector-valued) differential forms; in particular,
there is a naturally defined notion of vertical exterior derivative,
 \begin{equation} \label{eq:VEXTD1} % vertical exterior derivative
 \begin{array}{cccc}
  d_V^{} : & \Omega_V^{\,r}(P;\pi^* \hat{T})
           & \longrightarrow & \Omega_{~V}^{\,r+1}(P;\pi^* \hat{T}) \\[1mm]
           & \alpha & \longmapsto & d_V^{} \alpha
 \end{array}~,
\end{equation}
and of vertical Lie derivative along a vertical vector field $X$,
 \begin{equation} \label{eq:VLIED1} % vertical Lie derivative
 \begin{array}{cccc}
  \mathbb{L}_X^{} : & \Omega_V^{\,r}(P;\pi^* \hat{T}) & \longrightarrow
                    & \Omega_V^{\,r}(P;\pi^* \hat{T}) \\[1mm]
                    & \alpha & \longmapsto & \mathbb{L}_X^{} \alpha
 \end{array}~,
\end{equation}
which are defined by exactly the same formulas as in the standard case, namely
\begin{equation} \label{eq:VEXTD2} % vertical exterior derivative
 \begin{array}{rcl}
  (d_V^{} \alpha)(X_0^{},\ldots,X_r^{}) \!\!
  &=&\!\! {\displaystyle
           \sum_{i=0}^r \, (-1)^i \; X_i^{} \cdot
           \bigl( \alpha(X_0^{},\ldots,\hat{X}_i^{},\ldots,X_r^{}) \bigr)}
  \\[4mm]
  & & +\, {\displaystyle
           \sum_{0 \leqslant i < j \leqslant r} (-1)^{i+j} \,
           \alpha([X_i^{},X_j^{}],X_0^{},\ldots,
                  \hat{X}_i^{},\ldots,\hat{X}_j^{},\ldots,X_r^{})}
 \end{array}~,
\end{equation}
where $\; X_0^{},X_1^{},\ldots,X_r^{} \smin\, \mathfrak{X}_V^{}(P) \,$,
and
\begin{equation} \label{eq:VLIED2} % vertical Lie derivative
 (\mathbb{L}_X^{} \alpha)(X_1^{},\ldots,X_r^{})~
 =~X^{} \cdot \bigl( \alpha(X_1^{},\ldots,X_r^{}) \bigr) \, - \,
        \sum_{i=1}^r \, \alpha(X_1,\ldots,[X,X_i^{}],\ldots,X_r^{})~,
\end{equation}
where $\; X,X_1^{},\ldots,X_r^{} \smin\, \mathfrak{X}_V^{}(P) \,$:
this makes sense since $V\!P$ is an involutive distribution on~$P$.
Here and throughout the remainder of this section, the symbol $\cdot$
stands for the directional derivative of sections of~$\pi^*(\hat{T})$
along vertical vector fields: this makes sense since upon restriction
to each fiber, a vertical vector field is simply an ordinary vector
field on the fiber and a section of a vector bundle obtained as the
pull-back of a vector bundle over~$M$ becomes a function on the fiber
taking values in a fixed vector space.
\begin{dfn}~ \label{def:PLFB}
 A \textbf{polypresymplectic fiber bundle} is a fiber bundle\/~$P$ over an
 $n$-dimensional manifold\/~$M$ equipped with a vertical\/~$(k+1)$-form
 \[
  \hat{\omega} \in \Omega_{~V}^{\,k+1}(P;\pi^*(\hat{T}))
 \]
 of constant rank on the total space\/~$P$ taking values in (the pull-back 
 to\/~$P$ of) a fixed $\hat{n}$-dimensional vector bundle\/~$\hat{T}$ over 
 the same manifold\/~$M$, called the \textbf{polypresymplectic form along
 the fibers} of~$P$, or simply the \textbf{polypresymplectic form}, and
 said to be of \textbf{rank}\/~$N$, such that $\,\hat{\omega}$ is
 vertically closed,%
 \footnote{As in the symplectic case, the possible absence of the
 integrability condition $\, d_V^{} \hat{\omega} = 0 \,$ will be
 indicated by adding the term ``almost''.}
 \begin{equation}
  d_V^{} \hat{\omega}~=~0~,
 \end{equation}
 and such that at every point\/ $p$ of\/~$P$, $\hat{\omega}_p^{}$ is a
 polypresymplectic form of rank $N$ on the vertical space\/~$V_{\!p}^{} P$:
 this means that there exists a subspace\/ $L_p^{}$ of\/~$V_{\!p}^{} P$ of
 codimension\/~$N$, called the \textbf{polylagrangian subspace},%
 \footnote{The terminology, as well as the justification for using the
 definite article, stems from the fact that this subspace, if it exists,
 is more than just lagrangian (i.e., maximal isotropic) and that, as
 soon as either $\, \hat{n} > 1 \,$ or else $\, \hat{n} = 1 \,$ but
 then $\, N > k > 1$, or in other words, except when $\,\hat{\omega}$
 is an ordinary two-form or a volume form, it is necessarily unique.}
 \addtocounter{footnote}{-1}
 such that the ``musical map''
 \[
  \hat{\omega}_p^\flat :~V_{\!p}^{} P~\longrightarrow~
  \bwedge{k} \, V_p^* P \,\otimes\, \hat{T}_{\pi(p)}^{}
 \]
 given by contraction of~$\,\hat{\omega}_p^{}$ in its first argument,
 when restricted to~$L_p^{}$, yields a linear isomorphism
 \[
  L_p^{} \, / \, \mathrm{ker} \, \hat{\omega}_p^{}~\cong~
  \bwedge{k} L_p^\perp \,\otimes\, \hat{T}_{\pi(p)}^{}
 \]
 where\/ $L_p^\perp$ is the annihilator of\/~$L_p^{}$ in\/~$V_p^* P$.
 Moreover, it is assumed that the kernels\/~$\ker \, \hat{\omega}_p^{}$
 as well as the polylagrangian subspaces\/ $L_p^{}$ at the different
 points of\/~$P$ fit together \linebreak smoothly into distributions\/~%
 $\ker \, \hat{\omega}$ and\/ $L$ on\/~$P$: the latter is called the
 \textbf{polylagrangian distribution} of~$\,\hat{\omega}$. 
 If $\,\hat{\omega}$ is non-degenerate, we say that\/~$P$ is a
 \textbf{polysymplectic \linebreak fiber bundle} and $\,\hat{\omega}$
 is a  \textbf{polysymplectic form along the fibers} of\/~$P$, or
 simply a \textbf{polysymplectic form}.
 If\/ $M$ reduces to a point, we speak of a
 \textbf{poly(pre)symplectic manifold}.
 The case of main interest is when $\,\hat{\omega}$ is a $2$-form,
 i.e., $k=1$.
\end{dfn}
Thus the characteristic feature of a polysymplectic fiber bundle $P$
with a polysymplectic form $\,\hat{\omega}$ is the existence of a special 
subbundle~$L$ of its vertical bundle~$V\!P$ which is not only lagrangian 
(in particular, isotropic) but has the even stronger property that the 
``musical vector bundle homomorphism'' $\; \hat{\omega}^\flat :
V\!P \longrightarrow \bwedge{k} \, V^* P \,\otimes\, \pi^*(\hat{T}) \,$,
when restricted to~$L$, provides a vector bundle isomorphism
\begin{equation} \label{eq:MUSISO5}
 \hat{\omega}^\flat :~L~\stackrel{\cong}{\longrightarrow}~
 \bwedge{k} L^\perp \,\otimes\, \pi^*(\hat{T})~.
\end{equation}
As has been proved in Ref.~\cite{FG}, as soon as~$\, \hat{n} > 2$, $L$ is
necessarily involutive.

The following example provides what may be considered the ``standard model''
of a polysymplectic fiber bundle:
\begin{exe}~ \label{exe:MPFPL}
 Let $E$ be an arbitrary fiber bundle over an $n$-dimensional manifold~$M$,
 with projection $\, \pi_E : E \longrightarrow M$, and let $\hat{T}$ be a
 fixed $\hat{n}$-dimensional vector bundle over the same manifold~$M$.
 Consider the bundle
 \begin{equation} \label{eq:MPFPL}
  P~=~\bwedge{k} V^* E \,\otimes\, \pi_E^*(\hat{T})
 \end{equation}
 of vertical $k$-forms on~$E$ taking values in the pull-back of $\hat{T}$
 to~$E$, with projections \linebreak $\pi^k : P \longrightarrow E \,$ and
 $\, \pi = \pi_E \circ \pi^k : P \longrightarrow M$. Using the tangent
 map $\, T \pi^k : TP \longrightarrow TE$ \linebreak of~$\pi^k$ and its
 restriction $\, V \pi^k : VP \longrightarrow VE \,$ to the vertical
 bundles, we define the \linebreak \textbf{canonical $k$-form} on~$P$,
 which is a vertical $k$-form $\hat{\theta}$ on~$P$ taking values
 in~$\pi^*(\hat{T})$, by
 \begin{equation}
  \begin{array}{c}
   {\displaystyle
    \hat{\theta}_\alpha(v_1,\ldots,v_k)~
    =~\alpha \bigl( V_\alpha \pi^k \cdot v_1,\ldots,V_\alpha \pi^k \cdot v_k)}
   \\[1mm]
   \mbox{for $\, \alpha \in P \,$ and $\, v_1,\ldots,v_k \in V_\alpha P$}~.
  \end{array}
 \end{equation}
 Then $\; \hat{\omega} = - d_V \hat{\theta} \,$ is a polysymplectic
 $(k+1)$-form,  with polylagrangian distribution $\, L = \ker(T \pi^k)$
 (the vertical bundle for the projection to~$E$), contained in $\, VP =
 \ker(T \pi)$ (the vertical bundle for the projection to~$M$).
\end{exe}
When $\, k=1 \,$ and $\, \hat{n}=1$ (with the understanding that the
auxiliary vector bundle $\hat{T}$ is the trivial real line bundle
$\, M \times \mathbb{R}$), we have the ``standard model'' of a
symplectic fiber bundle.
In particular, when, in addition, $n=0$ (i.e., the base manifold $M$
is reduced to a single point), we recover the cotangent bundle of the
single fiber, which is an arbitrary manifold, as the ``standard model''
of a symplectic manifold.
On the other hand, when $\, k=1 \,$ and $\, \hat{n}=n-1$ (with the
understanding that the auxiliary vector bundle $\hat{T}$ is the
bundle $\bwedge{n-1} T^* M$ of $(n-1)$-forms on~$M$), $P$ can
be identified with the twisted dual $\vec{J}^{\oast} E$ of the
linearized jet bundle $\vec{J} E$ of~$E$ (which is the difference
vector bundle of the usual jet bundle $JE$ of~$E$), because
\[
 \vec{J} E~\cong~\pi_E^*(T^* M) \otimes VE
\]
implying
\[
 \vec{J}^\ast E~\cong~V^* E \otimes \pi_E^*(TM)
\]
for the common dual $\vec{J}^\ast E$ and
\[
 \vec{J}^{\oast} E~\cong~V^* E \otimes \pi_E^*(TM) \otimes
                       \pi_E^* \bigl( \bwedge{n} T^* M \bigr)
\]
for the twisted dual $\, \vec{J}^{\oast} E = \vec{J}^\ast E \otimes \pi_E^*
\bigl( \bwedge{n} T^* M \bigr) \,$, so we get a canonical isomorphism
\begin{equation}
 \vec{J}^{\oast} E~\cong~V^* E \otimes
                         \pi_E^* \bigl( \bwedge{n-1} T^* M \bigr)
\end{equation}
of vector bundles over~$E$.
This bundle plays an important role in the covariant hamiltonian formalism
of classical field theory~\cite{CCI,GIM,MG,FR}.

As a first application of the isomorphism~(\ref{eq:MUSISO5}) beyond those
discussed in Ref.~\cite{FG}, we show that, just as in the symplectic case,
it allows us to construct a polysymplectic version of the Bott connection. 
The idea is simple: start with the Bott connection in~$L^\perp$ as defined
in eqns~(\ref{eq:BOTTC1}) and~(\ref{eq:BOTTC2}) (with $\mathfrak{X}(M)$ 
replaced by $\mathfrak{X}_V^{}(P)$, $\Omega^1(M)$ replaced by $\Omega_V%
^{\,1}(P)$ and the common Lie derivative replaced by the vertical Lie 
derivative introduced at the beginning of this section) and take the
tensor product of its $k$-th exterior power with the trivial partial
linear connection in $\pi^*(\hat{T})$ along~$L$, to obtain a partial linear 
connection $\nabla^B$ in $\, \bwedge{k} L^\perp \otimes \pi^*(\hat{T}) \,$ 
along~$L$, which we call the \textbf{Bott connection in} $\, \bwedge{k} 
L^\perp \otimes\, \pi^*(\hat{T}) \,$: it is still given by a suitable 
restriction of the (vertical) Lie derivative of (vertical) forms;
explicitly,
\begin{equation} \label{eq:BOTTC5}
 \begin{array}{c}
  {\displaystyle
   (\nabla_{\!X}^B \alpha)(Y_1^{},\ldots,Y_k^{})~
   =~X \cdot \bigl( \alpha(Y_1^{},\ldots,Y_k^{}) \bigr) \, - \,
     \sum_{i=1}^k \, \alpha(Y_1^{},\ldots,[X,Y_i^{}],\ldots,Y_k^{})} \\[5mm]
  \mbox{for $\, X \in \Gamma(L)$, $\alpha \in
                      \Gamma(\bwedge{k} L^\perp \otimes \pi^*(\hat{T}))$,
            $Y_1^{},\ldots,Y_k^{} \in \mathfrak{X}_V^{}(P)$}~.
 \end{array}
\end{equation}
Now using the isomorphism~(\ref{eq:MUSISO5}), we can transfer it to a
partial linear connection $\nabla^B$ in~$L$ along~$L$ and arrive at
\begin{dfn}~ \label{dfn:PLBOTTC}
 Let\/ $P$ be an almost polysymplectic fiber bundle over a manifold\/~$M$
 with almost polysymplectic form $\,\hat{\omega}$ and involutive poly%
 lagrangian distribution\/~$L$. Then there exists a naturally defined
 partial linear connection\/ $\nabla^B$ in\/~$L$ along\/~$L$ which we
 call the \textbf{poly\-symplectic Bott connection}; explicitly, it is
 determined by the formula
 \begin{eqnarray} \label{eq:BCPL1}
  \begin{array}{c}
   {\displaystyle
    \hat{\omega} \bigl( \nabla_{\!X}^B Y , Z_1^{},\ldots,Z_k^{} \bigr)~
    =~X \cdot \bigl( \hat{\omega}(Y,Z_1^{},\ldots,Z_k^{}) \bigr) -
      \sum_{i=1}^k \, \hat{\omega}(Y,Z_1,\ldots,[X,Z_i^{}],\ldots,Z_k^{})}~
   \\[5mm]
   \mbox{for $\, X,Y \in \Gamma(L)$,
             $Z_1^{},\ldots,Z_k^{} \in \mathfrak{X}_V^{}(P)$}~.
  \end{array} \!\!\!
 \end{eqnarray}
\end{dfn}
As in the symplectic case, the polysymplectic Bott connection is flat, and
for its torsion we have the following analogue of Theorem~\ref{thm:TORBCSY}\,:
\begin{prp}~ \label{prp:TORBCPL}
 Let\/ $P$ be an almost polysymplectic fiber bundle over a manifold\/~$M$
 with almost polysymplectic form $\,\hat{\omega}$ and involutive poly%
 lagrangian distribution\/~$L$. Then if $\,\hat{\omega}$ is vertically
 closed, the polysymplectic Bott connection $\nabla^B$ has zero torsion.
 More generally, the torsion tensor\/~$T^B$ of\/~$\nabla^B$ is related
 to the vertical exterior derivative of~$\,\hat{\omega}$ by the formula
 \begin{equation} \label{eq:TBCPL1}
  \begin{array}{c}
   d_V^{} \hat{\omega} \, (X,Y,Z_1^{},\ldots,Z_k^{} \bigr)~
   =~\hat{\omega}(T^B(X,Y),Z_1^{},\ldots,Z_k^{})
   \\[2mm]
   \mbox{for $\, X,Y \in \Gamma(L)$,
             $Z_1^{},\ldots,Z_k^{} \in \mathfrak{X}_V^{}(P)$}~.
  \end{array}
 \end{equation}
 In particular, this implies that in a polysymplectic fiber bundle, the
 leaves of the poly\-lagrangian foliation are flat affine manifolds.
\end{prp}
\proof
 The only property that remains to be checked is eqn~(\ref{eq:TBCPL1}):
 this is a simple calculation using the fact that $\, \hat{\omega}(X,Y,
 \ldots) = 0 \,$ when $\, X,Y \in \Gamma(L) \,$ since $L$ is isotropic:
 \begin{eqnarray*}
 \lefteqn{\hat{\omega} \bigl( \nabla_{\!X}^B Y - \nabla_{\!Y}^B X - [X,Y] ,
                              Z_1^{},\ldots,Z_k^{} \bigr)} \hspace*{2cm} \\[5mm]
  &=&\!\! X \cdot \bigl( \hat{\omega}(Y,Z_1^{},\ldots,Z_k^{}) \bigr) - \,
          Y \cdot \bigl( \hat{\omega}(X,Z_1^{},\ldots,Z_k^{}) \bigr) \\
  & & \mbox{} - \, \sum_{i=1}^k \, (-1)^i \; Z_i^{} \cdot
      \bigl( \hat{\omega}(X,Y,Z_1^{},\ldots,\hat{Z}_i^{},\ldots,Z_k^{} \bigr)
  \\[2mm]
  & & \mbox{} - \, \hat{\omega}([X,Y],Z_1^{},\ldots,Z_k^{}) \\[2mm]
  & & \mbox{} - \, \sum_{i=1}^k \, (-1)^i \;
      \hat{\omega}([X,Z_i^{}],Y,Z_1,\ldots,\hat{Z}_i^{},\ldots,Z_k^{}) \\
  & & \mbox{} + \, \sum_{i=1}^k \, (-1)^i \;
      \hat{\omega}([Y,Z_i^{}],X,Z_1,\ldots,\hat{Z}_i^{},\ldots,Z_k^{}) \\
  & & \mbox{} + \, \sum_{1 \leqslant i < j \leqslant k}^{\phantom{k}} \,
      (-1)^{i+j} \; \hat{\omega}([Z_i^{},Z_j^{}],X,Y,Z_1^{},\ldots,
                                 \hat{Z}_i^{},\ldots,\hat{Z}_j^{},\ldots,Z_k^{})
  \\[4mm]
  &=&\!\! d_V^{} \hat{\omega} \bigl( X,Y,Z_1^{},\ldots,Z_k^{} \bigr)~.
 \end{eqnarray*}
\qed

Now we turn to polysymplectic connections. By analogy with the symplectic case,
the definition of the concept is more or less obvious, except that we cannot 
expect to obtain anything beyond partial linear connections along the vertical 
bundle.
\begin{dfn}~ \label{def:PLCON}
 A \textbf{poly(pre)symplectic connection} on an almost poly(pre)symplectic
 fiber bundle\/~$P$ over a manifold\/~$M$ with almost poly(pre)symplectic
 form $\,\hat{\omega}$ and poly\-lagrangian distribution\/~$L$ is a partial
 linear connection\/~$\nabla$ in the vertical bundle\/~$V\!P$ of\/~$P$
 along\/~$V\!P$ itself which preserves both $\,\hat{\omega}$ and\/~$L$;
 in particular, it satisfies\/ $\, \nabla \hat{\omega} = 0$, or explicitly,
 \begin{equation} \label{eq:PLCONN}
  \begin{array}{c}
   {\displaystyle
    X \cdot \bigl( \hat{\omega}(X_0^{},\ldots,X_k^{}) \bigr)~
    =~\sum_{i=0}^k \, \hat{\omega}(X_0^{},\ldots,\nabla_{\!X}^{} X_i^{},
                                         \ldots,X_k^{})} \\[5mm]
   \mbox{for $\, X,X_0^{},\ldots,X_k^{} \smin\, \mathfrak{X}_V^{}(P)$}~.
  \end{array}
 \end{equation}
\end{dfn}
Note that as soon as the polylagrangian distribution $L$ is unique, the
invariance of $\,\hat{\omega}$ under parallel transport with respect to~%
$\nabla$ already implies that of~$L$. As mentioned before,\footnotemark\,
the only exceptions to this situation can occur when $\,\hat{\omega}$ is
an ordinary two-form or a volume form: in these cases, invariance of~$L$
becomes a separate condition.

As in the symplectic case, the existence of torsion-free poly(pre)symplectic
connections imposes certain constraints.

\begin{prp}~\label{prp:PLCONN1}
 Let\/ $P$ be an almost poly(pre)symplectic fiber bundle over a manifold\/~$M$
 with almost poly(pre)symplectic form~$\,\hat{\omega}$ and polylagrangian
 distribution\/~$L$. Then if there exists a torsion-free poly(pre)symplectic
 connection\/ $\nabla$ on\/~$P$, $\,\hat{\omega}$ must be vertically closed
 and\/ $L$ must be involutive. More generally, the torsion tensor $T$ of a
 poly(pre)symplectic connection~$\nabla$ over $P$ is related to the 
 exterior derivative of $\,\hat{\omega}$ by the formula
 \begin{equation}
  d_V^{} \hat{\omega}(X_0^{},\ldots,X_r^{})~
  = \; - \sum_{0 \leqslant i < j \leqslant r}^{\phantom{r}} \, (-1)^{i+j} \;
         \hat{\omega}(T(X_i^{},X_j^{}),X_0^{},\ldots,\hat{X}_i^{},\ldots,
                                      \hat{X}_j^{},\ldots,X_r^{})~.
 \end{equation}
 Finally, if $\,\hat{\omega}$ is non-degenerate, the restriction of
 any torsion-free polysymplectic connection to\/~$L$ coincides with
 the polysymplectic Bott connection in\/~$L$.
\end{prp}
\begin{rmk}~
 As in the symplectic case, the last statement is valid under the same less
 restrictive assumptions on the torsion tensor~$T$ of~$\nabla$ as before:
 it suffices that $T(X,Y)$ should be along~$L$ as soon as $X$ or $Y$ is
 along~$L$.
\end{rmk}
\proof
 The first two statements follow directly from Lemma~\ref{lem:CONN2} 
 (eqn~(\ref{eq:LEMC2})) and Lemma~\ref{lem:CONN1} in Appendix~A.
 For the third statement, let us assume that $\,\hat{\omega}$ is non-%
 degenerate and $\nabla$ is any polysymplectic connection on~$P$ with
 torsion tensor~$T$ such  that $T(X,Y)$ is along~$L$ as soon as $X$
 or $Y$ is along~$L$. Then the claim is equivalent to the condition
 that for all $\, X,Y \in \Gamma(L) \,$ and $\, Z_1^{},\ldots,Z_k^{}
 \in \mathfrak{X}_V^{}(P)$,
 \[
  \hat{\omega}(\nabla_{\!X}^B Y,Z_1^{},\ldots,Z_k^{})~
  =~\hat{\omega}(\nabla_{\!X}^{} Y,Z_1^{},\ldots,Z_k^{})
 \]
 which can be derived by comparing eqn~(\ref{eq:BCPL1}) with the
 corresponding expression from eqn~(\ref{eq:PLCONN}) taking into
 account that, for $\, 1 \leqslant i \leqslant k$,
 \[
  \hat{\omega}(Y,Z_1^{},\ldots,[X,Z_i^{}],\ldots,Z_k^{})~
  =~\hat{\omega}(Y,Z_1^{},\ldots,\nabla_{\!X}^{} Z_i^{},\ldots,Z_k^{})
 \]
 since $L$ being isotropic and stable under $\nabla$, the expressions
 $\, \hat{\omega}(Y,Z_1^{},\ldots,\nabla_{\!Z_i}^{} X,\ldots,Z_k^{})$
 \linebreak
 and $\, \hat{\omega}(Y,Z_1^{},\ldots,T(X,Z_i^{}),\ldots, Z_k^{}) \,$
 vanish under these assumptions.
\qed

\noindent
Conversely, we can use a partition of unity argument to prove that the
conditions stated in Proposition~\ref{prp:PLCONN1} ($\hat{\omega}$ is
vertically closed and $L$ is involutive) are not only necessary but
also sufficient to guarantee existence of torsion-free polysymplectic
connections.
\begin{thm}~\label{thm:PLCONN1}
 Let\/~$P$ be a poly(pre)symplectic fiber bundle over a manifold\/~$M$ with
 poly(pre)symplectic form~$\,\hat{\omega}$ and involutive polylagrangian
 distribution\/~$L$. Then there exist\linebreak torsion-free poly(pre)%
 symplectic connections~$\,\nabla$ on\/~$P$, and the set of all such
 connections constitutes an affine space whose difference vector space
 can, for non-degenerate~$\,\hat{\omega}$, be identified with the space
 of all vertical tensor fields of rank\/~$k+2$ on\/~$P$ taking values in
 the auxiliary vector bundle\/~$\pi^*(\hat{T})$ which (a) have symmetry
 corresponding to the irreducible representation of the permutation
 group\/~$S_{k+2}$ given by the Young pattern
 \begin{center}
 \unitlength4mm
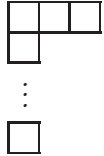
 \begin{picture}(7,5)
 \linethickness{0.2mm}
  \multiput(1.2,4)(1,0){3}{\framebox(0.95,0.95){}}
  \put(1.2,3){\framebox(0.95,0.95){}}
  \put(1.55,1.5){\vdots}
  \put(1.2,0){\framebox(0.95,0.95){}}
 \end{picture}
 \end{center}
 (with\/ $k$ boxes in the first column) and (b) vanish whenever at least
 two of their arguments are along\/~$L$.
\end{thm}
\proof
 Concerning existence, we can under the hypotheses of the theorem apply
 the polysymplectic Darboux theorem (see~\cite[Theorem~10]{FG}) to
 guarantee that locally (i.e., on a sufficiently small open neighborhood
 of each point of~$P$), there exists a torsion-free flat polysymplectic
 connection on~$P$: it is simply the partial linear connection in~$V\!P$
 along~$V\!P$ whose Christoffel symbols vanish identically in these
 coordinates. Now using a covering of~$P$ by such Darboux coordinate
 neighborhoods, passing to a locally finite refinement $(U_\alpha)%
 _{\alpha \in A}$, denoting the corresponding family of partial linear
 connections by $(\nabla_{\!\alpha})_{\alpha \in A}$ and choosing a
 partition of unity $(\chi_\alpha)_{\alpha \in A}$ subordinate to
 the open covering $(U_\alpha)_{\alpha \in A}$, we can define
 \[
  \nabla~=~\sum_{\alpha \in A} \, \chi_\alpha \, \nabla_{\!\alpha}~.
 \]
 Then it is clear that $\nabla$ is a partial linear connection in~$V\!P$
 along~$V\!P$, preserves~$\,\hat{\omega}$ as well as~$L$ and
 is torsion-free, since this is true for each $\nabla_{\!\alpha}$ and since
 the conditions of preserving a given differential form, of preserving a
 given vector subbundle and of being torsion-free are all local (i.e.,
 behave naturally under restriction to open subsets) as well as affine.%
 \footnote{With respect to curvature, the same comment as in Footnote~%
 \ref{fn:CURV} applies.}
 Regarding uniqueness, or rather the amount of non-uniqueness, we can write
 the difference between any partial linear connection $\nabla'$ in~$V\!P$
 along~$V\!P$ and a fixed torsion-free polysymplectic connection~$\nabla$
 on~$P$ in the form
 \[
  \nabla_{\!X}' Y~=~\nabla_{\!X}^{} Y \, + \, S(X,Y)~.
 \]
 Moreover, we introduce a (totally covariant) vertical tensor field
 $\,\hat{\omega}_S^{}$ of rank~$k+2$ on~$P$ taking values in $\hat{T}$
 which, for non-degenerate $\,\hat{\omega}$, carries exactly the same
 information as~$S$ itself, given by
 \[
  \hat{\omega}_S^{}(X,Y,Z_1^{},\ldots,Z_k^{})~
  =~\hat{\omega}(S(X,Y),Z_1^{},\ldots,Z_k^{})~.
 \]
 Obviously, $\hat{\omega}_S^{}$ is totally antisymmetric in its last $k$
 arguments and it is clear that $\nabla'$ will be torsion-free if and only
 if $S$ is symmetric, or equivalently, $\hat{\omega}_S^{}$ is symmetric in
 its first two arguments, that $\nabla'$ will preserve $\,\hat{\omega}$ if
 and only if $S$ satisfies the identity
 \[
  \hat{\omega}(S(X,Y),Z_1^{},\ldots,Z_k^{}) \, + \,
  \sum_{i=1}^k \, \hat{\omega}(Y,Z_1^{},\ldots,Z_{i-1}^{},S(X,Z_i^{}),
                               Z_{i+1}^{},\ldots,Z_k^{})~=~0~,
 \]
 or equivalently, $\hat{\omega}_S^{}$ satisfies the cyclic identity
 \[
  \hat{\omega}_S^{}(X,Y,Z_1^{},\ldots,Z_k^{}) \, - \,
  \sum_{i=1}^k \, \hat{\omega}_S^{}(X,Z_i^{},Z_1^{},\ldots,Z_{i-1}^{},Y,
                                    Z_{i+1}^{},\ldots,Z_k^{})~=~0~,
 \]
 and that $\nabla'$ will preserve~$L$ if and only if $S(X,Y)$ is along~$L$
 whenever $X$ or $Y$ are along~$L$, or equivalently, $\hat{\omega}_S^{}$
 vanishes whenever at least two of its arguments are along~$L$.
 Finally, it is well known that, together with symmetry in the
 first two arguments and antisymmetry in the last $k$ arguments,
 this cyclic identity identifies the tensor $\hat{\omega}_S^{}$
 as belonging to the irreducible representation of the permutation
 group $S_{k+2}$ given by the Young pattern stated in the theorem;
 see, e.g., \cite[p.~249]{Ha}.
\qed

\section{Multisymplectic connections}

We begin by stating the definition of a multisymplectic structure, as 
given in Ref.~\cite{FG}. To this end, we recall first of all that given 
a fiber bundle $P$ over a manifold~$M$ with bundle projection $\; \pi: P
\longrightarrow M$, an $r$-form on~$P$ is said to be $(r-s)$-horizontal, 
where $\, 0 \leqslant s \leqslant r$, if its contraction with more than 
$s$ vertical vectors vanishes; we shall in what follows denote the bundle 
of such forms by $\dbwedge{s}{r} \, T^* P$.
\begin{dfn}~ \label{def:MLFB}
 A \textbf{multipresymplectic fiber bundle} is a fiber bundle\/ $P$
 over an\linebreak $n$-dimensional manifold\/~$M$ equipped with a\/ 
 $(k+1-r)$-horizontal\/ $(k+1)$-form
 \[
  \omega \in \Gamma \bigl( \dbwedge{\;~r}{k+1} \, T^* P \bigr)
 \]
 of constant rank on the total space\/~$P$, where $\, 1 \leqslant r 
 \leqslant k+1 \,$ and $\, k+1-r \leqslant n$, called the 
 \textbf{multipresymplectic form} and said to be of \textbf{rank}\/ $N$ and 
 \textbf{horizontality degree}\/ $k+1-r$, such that $\omega$ is closed,%
 \footnote{Once again, the possible absence of the integrability condition
  $\, d \>\! \omega = 0 \,$ will be indicated by adding the term ``almost''.}
 \begin{equation}
  d \>\! \omega~=~0~.
 \end{equation}
 and such that at every point\/ $p$ of\/~$P$, $\omega_p^{}$ is
 a multipresymplectic form of rank\/~$N$ and horizontality degree\/~$k+1-r$
 on the tangent space\/~$T_{\!p}^{} P$: this means that there exists a sub%
 space\/ $L_p^{}$ of\/~$T_{\!p}^{} P$ contained in\/~$V_{\!p}^{} P$ and of
 codimension\/~$N$ there, called the \textbf{multilagrangian subspace},%
 \footnote{The terminology, as well as the justification for using the
 definite article, stems from the fact that this subspace, if it exists,
 is more than just lagrangian (in particular, isotropic) and, as soon
 as either $\, r < k+1 \,$ or else $\, r = k+1 \,$ and then $\, n=0$,
 $N > k > 1 \,$ and also either $\, k+1-r < n \,$ or else $\, k+1-r = n$
 \linebreak and then $\, N+n > k > n+1$, is necessarily unique.}
 \addtocounter{footnote}{-1}
 such that the ``musical map''
 \[
  \omega_p^\flat :~T_{\!p}^{} P~\longrightarrow~\bwedge{k} \, T_p^* P
 \]
 or
 \[
  \omega_p^\flat :~V_{\!p}^{} P~\longrightarrow~\dbwedge{\,r-1}{\;~k} \, T_p^* P
 \]
 given by contraction of~$\,\omega_p^{}$ in its first argument,
 when restricted to~$L_p^{}$, yields a linear isomorphism
 \[
  L_p^{} \, / \, \mathrm{ker} \, \omega_p^{}~\cong~
  \bwedge{k} L_p^\perp \,\cap\, \dbwedge{\,r-1}{\;~k} \, T_p^* P~.
 \]
 Moreover, it is assumed that the kernels\/~$\ker \, \omega_p^{}$ as well
 as the multilagrangian subspaces\/~$L_p^{}$ at the different points
 of\/~$P$ fit together smoothly into distributions\/ $\ker \, \omega$
 and\/ $L$ on\/~$P$: the latter is called the \textbf{multilagrangian
 distribution} of~$\,\omega$.
 If $\,\omega$ is non-degenerate, we say that\/ $P$ is a \textbf%
 {multisymplectic fiber bundle} and $\,\omega$ is a \textbf%
 {multisymplectic form}. \linebreak
 If\/ $M$ reduces to a point, we speak of a \textbf{multi(pre)symplectic
 manifold}.
 The case of main interest is when $\,\omega$ is an $(n-1)$-horizontal
 $(n+1)$-form, i.e., $\, k=n$, $r=2$.
\end{dfn}
Again, the characteristic feature of a multisymplectic fiber bundle~$P$
with a multisymplectic form $\,\omega$ is the existence of a special
subbundle~$L$ of its tangent bundle $TP$, contained in its
vertical bundle $VP$, which is not only lagrangian (in particular, iso%
tropic) but has the even stronger property that the ``musical vector bundle
homomorphism'' $\; \omega^\flat : TP \longrightarrow \bwedge{k} \, T^* P \;$
or $\; \omega^\flat : VP \longrightarrow \dbwedge{\,r-1}{\,~k} \, T^* P \,$,
when restricted to~$L$, provides a vector bundle isomorphism
\begin{equation} \label{eq:MUSISO6}
 \omega^\flat :~L~\stackrel{\cong}{\longrightarrow}~
 \bwedge{k} L^\perp \,\cap\, \dbwedge{\,r-1}{\;~k} \, T^* P~.
\end{equation}
As has been proved in Ref.~\cite{FG}, as soon as~$\, \binom{n}{k+1-r} > 2$,
$L$ is necessarily involutive.

Again, the following example provides a ``standard model'' of a
multisymplectic fiber bundle:
\begin{exe}~ \label{exe:MPFML}
 Let $E$ be an arbitrary fiber bundle over an $n$-dimensional manifold~$M$,
 with projection $\, \pi_E : E \longrightarrow M$. Consider the bundle
 \begin{equation} \label{eq:MPFML}
  P~=~\dbwedge{\,r-1}{\;~k} T^* E
 \end{equation}
 of $(k+1-r)$-horizontal $k$-forms on~$E$, where $\, 1 \leqslant r \leqslant
 k+1 \,$ and $\, k+1-r \leqslant n$, with projections $\, \pi_{r-1}^{~k}:
 P \longrightarrow E \,$ and $\, \pi = \pi_E \circ \pi_{r-1}^{~k}: P
 \longrightarrow M$. Using the tangent map $\, T \pi_{r-1}^{~k} : TP
 \longrightarrow TE \,$ of~$\pi_{r-1}^{~k}$, we define the \textbf%
 {canonical $k$-form} on~$P$, which is a $(k+1-r)$-horizontal $k$-form
 $\theta$ on~$P$, by
 \begin{equation}
  \begin{array}{c}
   {\displaystyle
    \theta_\alpha^{}(v_1^{},\ldots,v_k^{})~
    =~\alpha \bigl( T_\alpha^{} \pi_{r-1}^{~k} \cdot v_1^{},\ldots,
                    T_\alpha^{} \pi_{r-1}^{~k} \cdot v_k^{})} \\[1mm]
   \mbox{for $\, \alpha \in P \,$ and
              $\, v_1^{},\ldots,v_k^{} \in T_\alpha^{} P$}~.
  \end{array}
 \end{equation}
 Then $\; \omega = - d{\theta} \,$ is a multisymplectic $(k+1)$-form,
 with multilagrangian distribution $\, L = \ker(T \pi_{r-1}^{~k})$
 (the vertical bundle for the projection to~$E$), contained in
 $\, VP = \ker(T \pi)$ (the vertical bundle for the projection
 to~$M$).
\end{exe}
When $\, k=n \,$ and $\, r=2$, we have the ``standard model'' of a
multisymplectic fiber bundle since in this case, as is explicitly
demonstrated in the literature~\cite{CCI,GIM,MG}, $P$ can be
identified with the twisted dual $J^{\ostar} E$ of the jet bundle
$JE$ of~$E$, i.e., we have a canonical isomorphism
\begin{equation}
 J^{\ostar} E~\cong~\dbwedge{\,1}{\,n} \, T^* E
\end{equation}
of vector bundles over~$E$.
This bundle plays a central role in the covariant hamiltonian formalism
of classical field theory~\cite{CCI,GIM,MG,FR}.

Again, as a first application of the isomorphism~(\ref{eq:MUSISO6}) beyond
those discussed in Ref.~\cite{FG}, we show that it allows us to construct a
multisymplectic version of the Bott connection. Namely, consider the $k$-th 
exterior power of the Bott connection in~$L^\perp$, as defined in eqns~%
(\ref{eq:BOTTC1}) and~(\ref{eq:BOTTC2}) (with $\mathfrak{X}(M)$ replaced
by $\mathfrak{X}(P)$), which is a partial linear connection $\nabla^B$ in 
$\bwedge{k} L^\perp$ along~$L$: it is still given by a suitable restriction
of the Lie derivative of forms; explicitly,
\begin{equation} \label{eq:BOTTC7}
 \begin{array}{c}
  {\displaystyle
   (\nabla_{\!X}^B \alpha)(Y_1^{},\ldots,Y_k^{})~
   =~X \cdot \bigl( \alpha(Y_1^{},\ldots,Y_k^{}) \bigr) \, - \,
     \sum_{i=1}^k \, \alpha(Y_1^{},\ldots,[X,Y_i^{}],\ldots,Y_k^{})} \\[5mm]
  \mbox{for $\, X \in \Gamma(L)$, $\alpha \in \Gamma(\bwedge{k} L^\perp)$,
            $Y_1^{},\ldots,Y_k^{} \in \mathfrak{X}(P)$}~.
 \end{array}
\end{equation}
Now noting that it preserves $\, \bwedge{k} L^\perp \,\cap\,
\dbwedge{\,r-1}{\;~k} \, T^* P \,$ (the expression in eqn~%
(\ref{eq:BOTTC7}) vanishes if at least $r$ of the vector fields
$\, Y_1^{},\ldots,Y_k^{} \,$ are vertical, since $\, L \subset VP \,$
and $VP$ is involutive), we can use the isomorphism~(\ref{eq:MUSISO6})
to transfer it to a partial linear connection $\nabla^B$ in~$L$ along~$L$
and arrive at
\begin{dfn}~ \label{dfn:MLBOTTC}
 Let\/~$P$ be an almost multisymplectic fiber bundle over a manifold\/ $M$
 with almost multisymplectic form $\,\omega$ and involutive multilagrangian
 distribution\/~$L$. Then there exists a naturally defined partial linear
 connection\/ $\nabla^B$ in\/~$L$ along\/~$L$ which we call the \textbf%
 {multisymplectic Bott connection}; explicitly, it is determined by
 the formula
 \begin{eqnarray} \label{eq:BCML1}
  \begin{array}{c}
   {\displaystyle
    \omega \bigl( \nabla_{\!X}^B Y , Z_1^{},\ldots,Z_k^{} \bigr)~
    =~X \cdot \bigl( \omega(Y,Z_1^{},\ldots,Z_k^{}) \bigr) -
      \sum_{i=1}^k \, \omega(Y,Z_1,\ldots,[X,Z_i^{}],\ldots,Z_k^{})}~
   \\[5mm]
   \mbox{for $\, X,Y \smin\, \Gamma(L)$,
             $Z_1^{},\ldots,Z_k^{} \smin\, \mathfrak{X}(P)$}~.
  \end{array} \!\!\!
 \end{eqnarray}
\end{dfn}
As in the symplectic case, the multisymplectic Bott connection is flat,
and for its torsion we have the following analogue of Theorem~%
\ref{thm:TORBCSY}.
\begin{prp}~ \label{prp:TORBCML}
 Let\/ $P$ be an almost multisymplectic fiber bundle over a manifold\/~$M$
 with almost multisymplectic form $\,\omega$ and involutive multilagrangian
 distribution\/~$L$. Then if\/$\,\omega$ is closed, the multisymplectic
 Bott connection $\nabla^B$ has zero torsion. More generally, the torsion
 tensor\/ $T^B$ of\/~$\nabla^B$ is related to the exterior derivative 
 of\/~$\,\omega$ by the formula
  \begin{equation} \label{eq:TBCML1}
  \begin{array}{c}
   d \>\! \omega \, (X,Y,Z_1^{},\ldots,Z_k^{} \bigr)~
   =~\omega(T^B(X,Y),Z_1^{},\ldots,Z_k^{})
   \\[2mm]
   \mbox{for $\, X,Y \in \Gamma(L)$,
             $Z_1^{},\ldots,Z_k^{} \in \mathfrak{X}(P)$}~.
  \end{array}
 \end{equation}
 In particular, this implies that in a multisymplectic fiber bundle,
 the leaves of the multi\-lagrangian foliation are flat affine manifolds.
\end{prp}
\proof
 The proof is identical to the proof of the Proposition~\ref{prp:TORBCPL} and 
 will therefore not be repeated here.
\qed

Now we turn to multisymplectic connections, whose definition is analogous
to the ones given previously, the main difference being that these are full 
connections and not just partial ones.
\begin{dfn}~ \label{def:MLCON}
 A \textbf{multi(pre)symplectic connection} on an almost multi\-(pre)\-%
 symplectic fiber bundle\/ $P$ over a manifold\/~$M$ with almost multi%
 (pre)symplectic form $\,\omega$ and multilagrangian distribution\/~$L$
 is a linear connection\/ $\nabla$ in~$P$ which preserves both $\,\omega$
 and~$L$, as well as the vertical bundle\/~$V\!P$ of~$P$; in particular,
 it satisfies\/ $\, \nabla \omega = 0$, or explicitly,
 \begin{equation} \label{eq:MLCONN}
  \begin{array}{c}
   {\displaystyle
    X \cdot \bigl( \omega(X_0^{},\ldots,X_k^{}) \bigr)~
    =~\sum_{i=0}^k \, \omega(X_0^{},\ldots,\nabla_{\!X}^{} X_i^{},
                                         \ldots,X_k^{})} \\[5mm]
   \mbox{for $\, X,X_0^{},\ldots,X_k^{} \smin\, \mathfrak{X}(P)$}~.
  \end{array}
 \end{equation}
\end{dfn}
From the point of view of fiber bundle theory, the requirement that the
vertical bundle should be invariant under parallel transport with respect
to~$\nabla$ is a natural consistency condition: it is necessary in order
that parallel transport maps points in the same fiber to points in the
same fiber. Concerning invariance of the multilagrangian distribution,
we note as before that as soon as $L$ is unique, the invariance of
$\,\omega$ under parallel transport with respect to~$\nabla$ already
implies that of~$L$. In the few exceptional cases where this uniqueness
does not prevail,\footnotemark\ invariance of~$L$ becomes a separate
condition.

Finally, we note that the existence of torsion-free multisymplectic
connections imposes the same kind of constraints as before and that
when these constraints are satisfied, such connections can be completely
classified. We just give the statements and omit the proofs since these
are obtained by almost literally repeating those of Proposition~%
\ref{prp:PLCONN1} and Theorem~\ref{thm:PLCONN1} above.
\begin{prp}~ \label{prp:MLCONN1}
 Let\/ $P$ be an almost multi(pre)symplectic fiber bundle over a manifold\/~%
 $M$ with almost multi(pre)symplectic form~$\,\omega$ and multilagrangian
 distribution\/~$L$. Then if there exists a torsion-free multi(pre)symplectic
 connection\/ $\nabla$ on\/~$P$, $\omega$ must be closed and\/ $L$ must be
 involutive.
 More generally, the torsion tensor\/~$T$ of a multi(pre)symplectic
 connection\/~$\nabla$ on\/~$P$ is related to the exterior derivative
 of~$\,\omega$ by the formula
 \begin{eqnarray}
  \begin{array}{c}
   {\displaystyle
    d\omega(X_0^{},\ldots,X_k^{})~
    = \; - \sum_{0 \leqslant i < j \leqslant k} \, (-1)^{i+j} \;
      \omega(T(X_i^{},X_j^{}),X_0^{},\ldots,\hat{X}_i^{},\ldots,
                              \hat{X}_j^{},\ldots,X_k^{})} \\[5mm]
   \mbox{for $\, X_0^{},\ldots,X_k^{} \in \mathfrak{X}(P)$}~.
  \end{array}
 \end{eqnarray}
 Finally, if $\,\omega$ is non-degenerate, the restriction of any
 torsion-free multisymplectic connection to\/~$L$ coincides with
 the multisymplectic Bott connection in\/~$L$.
\end{prp}
\proof
 Analogous to that of Proposition~\ref{prp:PLCONN1}.
\qed
\begin{thm}~ \label{thm:MLCONN1}
 Let\/ $P$ be a multi(pre)symplectic fiber bundle over a manifold\/~$M$
 with multi(pre)symplectic form~$\,\omega$ and involutive multilagrangian
 distribution\/~$L$. Then there exist torsion-free multi(pre)symplectic
 connections~$\,\nabla$ on\/~$P$, and the set of all such connections
 constitutes an affine space whose difference vector space can, for
 non-degenerate~$\,\omega$, be identified with the space of all tensor
 fields of rank\/~$k+2$ on\/~$P$ which (a) have symmetry corresponding
 to the irreducible representation of the permutation group\/~$S_{k+2}$
 given by the Young pattern
 \begin{center}
 \unitlength4mm
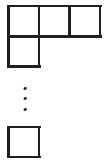
 \begin{picture}(7,5)
 \linethickness{0.2mm}
  \multiput(1.2,4)(1,0){3}{\framebox(0.95,0.95){}}
  \put(1.2,3){\framebox(0.95,0.95){}}
  \put(1.55,1.5){\vdots}
  \put(1.2,0){\framebox(0.95,0.95){}}
 \end{picture}
 \end{center}
 (with\/ $k$ boxes in the first column) and (b) vanish whenever at least
 $r+1$ of their arguments are vertical or at least two of their arguments
 are along\/~$L$.
\end{thm}
\proof
 Analogous to that of Theorem~\ref{thm:PLCONN1}.
\qed

\section{Structure Theorems}

In the previous two sections, we have presented `` standard models'' for
polysymplectic and multisymplectic fiber bundles: they are certain bundles
of forms built over a given fiber bundle, much in the same way as the
cotangent bundle of a given manifold is the ``standard model'' of a
symplectic manifold.
But of course we may wonder whether there are other interesting examples,
based on other methods.
In particular, a natural question to ask is whether there exists a
polysymplectic or multisymplectic analogue not only of the cotangent
bundle construction, but also of the coadjoint orbit construction of
symplectic geometry.

An alternative approach consists in looking at the converse question,
which in the context of symplectic geometry can be stated as follows:
How can we characterize, among all symplectic manifolds, those which
(up to a symplectomorphism) are cotangent bundles?
As it turns out, this issue is solved by Weinstein's
tubular neighborhood theorem.

As an initial step, we mention two conditions that are obviously necessary:
for a symplectic manifold $P$ to be symplectomorphic to the cotangent bundle
$T^* Q$ of some other manifold~$Q$, it must be exact (the de Rham cohomology
class of its symplectic form must vanish), and it must admit a lagrangian
foliation, whose leaves are of course the cotangent spaces.
Note that each of these conditions already excludes most of the
interesting coadjoint orbits, such as Souriau's $2$-sphere and,
more generally, all (co)adjoint orbits of compact semisimple
Lie groups, which are also K\"ahler manifolds.
But there are at least two other aspects that
turn out to be important.

The first aspect is that $P$ admits not only a lagrangian foliation
but also lots of submanifolds complementary to it: these submanifolds,
which may or may not be lagrangian, are the graphs of $1$-forms.%
\footnote{We use the term ``complementary'' as a stronger version of
the term ``transversal'': given two sub\-mani\-folds $X_1$ and $X_2$
of a manifold~$X$ and a point $x$ in their intersection, we say that
they are transversal at~$x$ if $\, T_x X = T_x X_1 + T_x X_2 \,$
and are complementary at~$x$ if $\, T_x X_1 \oplus T_x X_2 = T_x X$.
In the literature, a submanifold complementary to the leaves of a
foliation is often called a ``cross section'' of that foliation~--
a term inspired by fiber bundle theory when the submanifold is the
graph of some map.}
(As is well known, such a graph is a lagrangian submanifold if and
only if the corresponding $1$-form is closed.)
Note, however, that even though there are many such complementary
submanifolds (there are even many of them passing through each point
of~$P$), they are isolated, i.e., there is no canonical way to make
them come in families that would form a second foliation complementary
to the first one.
But at any rate, they are natural candidates for a manifold~$Q$
satisfying $\, P \cong T^* Q$.

The second aspect is that the lagrangian foliation is not arbitrary
but is simple, i.e., the quotient space of leaves can be given the
structure of a manifold such that the canonical projection becomes
a surjective submersion.
Again, this quotient space is a natural candidate for a manifold~$Q$
satisfying $\, P \cong T^* Q$.

Weinstein's tubular neighborhood theorem deals with the converse question:
Suppose that $P$ is a symplectic manifold, with symplectic form~$\,\omega$,
which admits a simple lagrangian foliation $\mathcal{F}$ (i.e., a lagrangian
foliation whose leaves are the connected components of the level sets of a
surjective submersion from~$P$ onto some other manifold), and let $Q$ be
any submanifold of~$P$ complementary to~$\mathcal{F}$.%
\footnote{Note that the condition that $Q$ should be complementary to~%
$\mathcal{F}$ is only local: it does not guarantee that the intersection
of~$Q$ with every leaf of~$\mathcal{F}$ reduces to a single point.
All it implies is that this intersection must be discrete and hence
at most countable, but it can contain many distinct points or, at
the other extreme, even be empty.}
In its original version~\cite{We}, the theorem states that if $Q$ is
lagrangian, then there is a tubular neighborhood of $Q$ in $P$ which is
symplectomorphic to a neighborhood of the zero section of the cotangent
bundle $T^*Q$ of~$Q$. This result is easily generalized to the case when
$Q$ is not lagrangian: it is enough to substitute the standard symplectic
form $-d\theta$ on~$T^* Q$ by a modified symplectic form $\, -d\theta +
\tau^* \omega_Q \,$ where $\tau$ is the canonical projection of~$T^* Q$
to~$Q$ and $\omega_Q$ is the restriction of~$\omega$ to~$Q$; see~\cite{DR}.
A global version of this result was given by Thompson~\cite{Th}, under the
hypothesis that the leaves of~$\mathcal{F}$ are simply connected and
geodesically complete: in this case, the manifold~$P$ becomes an affine
fiber bundle over the quotient manifold~$P/\mathcal{F}$ whose difference
vector bundle is its cotangent bundle $T^*(P/\mathcal{F})$, and therefore
there exist submanifolds $Q$ of~$P$ complementary to~$\mathcal{F}$ which
satisfy $\, Q \cong P/\mathcal{F}$, i.e., which meet every leaf of~%
$\mathcal{F}$ in precisely one point.%
\footnote{This follows from the fact that an affine fiber bundle always
admits global sections (which is easy to prove using partitions of unity).}
Furthermore, when we choose one such submanifold~$Q$, we get a global
symplectomorphism from~$P$ onto $T^* Q$ that takes $\omega$ to~$-d\theta$
or, more generally, to~$\, -d\theta + \tau^* \omega_Q \,$, as before.
The problem with the approach of~\cite{Th} is that it is not intrinsic,
since the structure of~$P$ as an affine bundle over~$P/\mathcal{F}$ and
the meaning of geodesic completeness of the leaves seem to depend on the
choice of additional ingredients (the author uses an auxiliary riemannian
metric, or rather its Levi-Civita connection).
 
In the remainder of the paper, we shall not only give a much more
transparent proof of all these theorems, but we shall also show that
this allows us to generalize them, without any additional effort, to
the setting of polysymplectic and multisymplectic geometry, where they
become natural structure theorems since these geometries come with an
in-built lagrangian foliation, right from the start.
The main natural ingredients used in our proofs, whose importance
in this context seems to have been underestimated in the past, are
(a) the Bott connection and (b) the concept of Euler vector field.

\section{Simple foliations by flat affine manifolds}

The main technical tool, which in what follows will be employed in
various different contexts and which therefore deserves to be treated
separately, in order to avoid unnecessary repetitions, is the notion
of a simple foliation of a manifold by flat affine submanifolds.

Initially, suppose that $P$ is any manifold. According to the Frobenius
theorem, a foliation $\mathcal{F}$ of~$P$ corresponds to an involutive
distribution $L$ on~$P$, such that for every point $p$ in~$P$,
\begin{equation} \label{eq:DIS}
 L_p~=~T_p \, \mathcal{F}_p~,
\end{equation}
where $\mathcal{F}_p$ is the leaf of~$\mathcal{F}$ passing through $p$.
Such a foliation is called \emph{simple}\/ if its leaves are the
connected components of the level sets of a surjective submersion
$\, \pi: P \longrightarrow \bar{P}$, that is, for $\, p \in P \,$
and $\, \bar{p} \in \bar{P} \,$ with $\, \bar{p} = \pi(p)$, we have
\begin{equation}
 \mathcal{F}_p~
 =~\mbox{connected component of $\pi^{-1}(\bar{p})$ containing $p$}~.
\end{equation}
In particular, a fiber bundle with connected fibers is a simple foliation.
It is also obvious that the leaves of a simple foliation are closed embedded
(and not just immersed) sub\-mani\-folds. Finally, given any surjective
submersion $\, \pi: P \longrightarrow \bar{P}$, we can always decompose
the projection $\pi$ into the composition of two projections,
\begin{equation} \label{eq:PRFOLH1}
 P~\longrightarrow~P/\mathcal{F}~\longrightarrow~\bar{P}~,
\end{equation}
where the first is a surjective submersion with connected fibers and the
second is a local diffeomorphism.

Given an arbitrary surjective submersion $\, \pi: P \longrightarrow \bar{P}$,
there are two special types of vector fields on the manifold~$P$: vertical
vector fields  and, more generally, projectable vector fields:
\begin{dfn}~ \label{def:CVPROJ}
 Let\/ $P$ and\/ $\bar{P}$ be manifolds and $\, \pi: P \longrightarrow
 \bar{P} \,$ be a surjective submersion. A vector field\/ $X$ on\/~$P$
 is said to be \textbf{vertical} (with respect to\/~$\pi$) if for any
 $\, p \in P$, $T_p \pi \cdot X(p) = 0$, and is said to be \textbf%
 {projectable} (with respect to\/~$\pi$) if for any $\, p_1,p_2 \in P \,$
 with $\, \pi(p_1) = \pi(p_2)$, $T_{p_1^{}} \pi \cdot X(p_1^{}) =
 T_{p_2^{}} \pi \cdot X(p_2^{})$.
\end{dfn}
If $X$ is a projectable vector field on~$P$, then it is clear
that for any $\, \bar{p} \in \bar{P}$, there exists a unique
vector $\, \bar{X}(\bar{p}) \in T_{\bar{p}} \bar{P} \,$ such
that $\, T_p \pi \cdot X(p) = \bar{X}(\bar{p}) \,$ for all
$\, p \in P \,$ with $\, \pi(p) = \bar{p}$, and using local
charts for~$P$ and~$\bar{P}$ in which the submersion $\pi$
is represented by a constant projection, we can check that
since $X$ is smooth, so is $\bar{X}$. Thus we can characterize
a projectable vector field as a vector field $X$ on~$P$ which
can be pushed forward by~$\pi$ to a (unique) vector field
$\bar{X}$ on~$\bar{P}$, to which it is $\pi$-related,%
\footnote{Recall that given any smooth map $\, f: M \longrightarrow N \,$
between manifolds $M$ and $N$, two vector fields $X$ on~$M$ and $Y$ on~$N$
are said to be $f$-related if for every point $m$ of~$M$, we have
$\, T_m f \cdot X(m) = Y(f(m))$. An elementary but important theorem,
used constantly, states that if $X_1$ on~$M$ and $Y_1$ on~$N$ are
$f$-related and $X_2$ on~$M$ and $Y_2$ on~$N$ are also $f$-related,
then their Lie brackets, $[X_1,X_2]$ on~$M$ and $[Y_1,Y_2]$ on~$N$,
are $f$-related.}
and a vertical vector field as a projectable vector field which,
when pushed forward by~$\pi$, gives zero. This implies immediately
that in the Lie algebra $\mathfrak{X}(P)$ of vector fields on~$P$,
the projectable vector fields form a Lie subalgebra $\mathfrak{X}_P(P)$
and the vertical vector fields form an ideal $\mathfrak{X}_V(P)$ within
this Lie subalgebra, i.e.
\begin{equation} \label{eq:CVPROJ1}
 \mbox{$Y,Z$ projectable}~~\Longrightarrow~~
 \mbox{$[Y,Z]$ projectable}~,
\end{equation}
\begin{equation} \label{eq:CVPROJ2}
 \mbox{$X$ vertical, $Y$ projectable}~~\Longrightarrow~~
 \mbox{$[X,Y]$ vertical}~.
\end{equation}
Finally, we have%
\begin{lem}~ \label{lem:CVPROJ1}
 Let\/ $P$ and\/ $\bar{P}$ be manifolds and $\, \pi: P \longrightarrow
 \bar{P} \,$ be a surjective submersion. Then every vector field\/
 $\bar{X}$ on\/~$\bar{P}$ is the push-forward of some projectable
 vector field\/ $X$ on\/~$P$ by\/~$\pi$.
\end{lem}
\proof
 Using local charts of~$P$ and~$\bar{P}$ where the submersion $\pi$ is
 represented by a constant projection, we see that every point $p$ of~$P$
 has an open neighborhood $U_p$ on which we can construct a vector field
 $X_p$ that projects to~$\bar{X}\big|_{\pi(U_p)}$. Choosing a locally
 finite refinement $(U_i)_{i \in I}$ of the open covering $(U_p)_{p \in P}$
 of~$P$ and a subordinate partition of unity $(\chi_i)_{i \in I}$, we can
 define a vector field $X$ on~$P$ by
 \[
  X~=~\sum_{i \in I} \, \chi_i \, X_{p(i)}
 \]
 and verify that its projects to~$\bar{X}$.
\qed

\pagebreak

\noindent
Thus we obtain the following exact sequence of Lie algebras:
\begin{equation}
 0~\longrightarrow~\mathfrak{X}_V(P)~\longrightarrow~\mathfrak{X}_P(P)~
   \longrightarrow~\mathfrak{X}(\bar{P})~\longrightarrow~0~.
\end{equation}
For later use, we note the following corollary:
\begin{lem}~ \label{lem:CVPROJ2}
 Let\/ $P$ and \/ $\bar{P}$ be manifolds and $\, \pi: P \longrightarrow
 \bar{P} \,$ be a surjective submersion. Then for every tangent vector
 $\, u \in T_p P$, there is a projectable vector field\/ $X$ on\/~$P$
 such that $\, X(p) = u$.
\end{lem}

Now we turn to the main subject of this section: the study of manifolds $P$
equipped with a simple foliation~$\mathcal{F}$ whose leaves are flat affine
submanifolds of~$P$. This means that the involutive distribution~$L$ tangent
to~$\mathcal{F}$, according to eqn~(\ref{eq:DIS}), is endowed with a partial
linear connection $\nabla$ with vanishing curvature and torsion. In this case,
we can define two special types of fields on~$P$, both of which are vertical
(i.e., along~$\mathcal{F}$, or~$L$) that play an important role: the covariant
constant vector fields and the Euler vector fields:
\begin{dfn}~ \label{def:EULVF}
 Let\/ $P$ be a manifold equipped with a simple foliation\/ $\mathcal{F}$
 with involutive tangent distribution\/~$L$, and let\/ $\nabla$ be a partial
 linear connection in\/~$L$ along\/~$L$ with vanishing curvature and torsion.
 We say that a vector field\/~$X$ tangent to\/~$\mathcal{F}$ is \textbf%
 {covariantly constant} along the leaves of\/~$\mathcal{F}$, or simply
 covariantly constant, if for any vector field\/~$Z$ tangent to\/~%
 $\mathcal{F}$, we have
 \begin{equation} \label{eq:CVCVF}
  \nabla_{\!Z}^{} X~=~0~.
 \end{equation}
 We also say that a vector field\/~$\Sigma$ tangent to\/~$\mathcal{F}$
 is an \textbf{Euler vector field} if for any vector field\/~$Z$ tangent
 to\/~$\mathcal{F}$, we have
 \begin{equation} \label{eq:EULVF}
  \nabla_{\!Z}^{} \Sigma~=~Z~.
 \end{equation}
\end{dfn}
The standard situation where these types of vector fields can be defined
naturally is on the total space of a vector bundle: in this case, there is
a preferred Euler vector field, namely, the one that vanishes on the zero
section.
However, the same construction also works for affine bundles~-- although
in this case, we lose uniqueness of the Euler vector field, since the
notion of the zero section has disappeared. In general, Definition~%
\ref{def:EULVF} implies immediately that the sum of a covariantly
constant vector field and an Euler vector field is an Euler vector
field, and conversely, the difference between two Euler vector fields
is a covariantly constant vector field, so the Euler vector fields
constitute an affine space whose difference vector space is the space
of covariantly constant vector fields. It is also clear that both types
of vector fields are uniquely determined by their value at a single
point of each leaf, and using local coordinate systems adapted to
the surjective submersion $\, P \longrightarrow P/\mathcal{F} \,$
in which the Christoffel symbols of the connection $\nabla$ vanish
identically, we can prove that both always exist, at least locally.

Completing the ``menu'' of ingredients, suppose now that $Q$ is a
submanifold of~$P$ complementary to~$\mathcal{F}$, i.e., for every
point $q$ of~$Q$, we have
\begin{equation} \label{eq:SVCOMP1}
 T_q P~=~T_q Q \,\oplus\, T_q \mathcal{F}_q~=~T_q Q \,\oplus\, L_q~.
\end{equation}
Note that this condition of complementarity does not necessarily imply
that $Q$ must intersect all leaves. However, considering again the
surjective submersion $\, \pi: P \longrightarrow P/\mathcal{F}$,
it does imply that every point of~$Q$ has an open neighborhood
in~$P$ whose intersection with~$Q$ is a local section of~$\pi$ and,
hence, that $\pi(Q)$ is open in~$P/\mathcal{F}$. Moreover, it also
implies that the inclusion of~$Q$ in~$P$, followed by the projection
$\pi$, as a map
\begin{equation} \label{eq:PRFOLH2}
 Q~\longrightarrow~P/\mathcal{F}
\end{equation}
is a local diffeomorphism onto its image, which is an open submanifold
of~$P/\mathcal{F}$. Thus, replacing $P$ by its open submanifold
$\pi^{-1}(\pi(Q))$ and $\mathcal{F}$ by restriction to this submanifold,
we can assume without loss of generality that $Q$ intersects all leaves,
i.e., that the map~(\ref{eq:PRFOLH2}) is surjective.

With these preliminaries out of the way, we want to show how to build,
using the geodesic flow with respect to the connection~$\nabla$ that
radially emanates from~$Q$, a canonical local diffeomorphism, denoted
by $\exp_Q$ and adequately called the \emph{exponential}, between
the vector bundle $L\big|_Q$ and the manifold~$P$. More precisely,
if for $\, q \in Q \,$ and $\, u_q \in L_q$, the geodesic in~%
$\mathcal{F}_q$ with initial position~$q$ and initial velocity~$u_q$
is (momentarily) denoted by $F(\,.\,;u_q)$, the map
\begin{equation} \label{eq:EXPON1}
 \exp_Q:~\mathrm{Dom}(\exp_Q)~\longrightarrow~P
\end{equation}
with domain given by
\begin{equation} \label{eq:EXPON2}
 \mathrm{Dom}(\exp_Q)~=~\bigcup_{q \in Q}~\{ \, u_q \in L_q \; | \;
                                             F(1;u_q)\ \text{exists} \, \}
\end{equation}
is defined by
\begin{equation} \label{eq:EXPON3}
 \exp_Q(u_q)~=~F(1;u_q)~.
\end{equation}
This allows us to immediately get rid of the symbol $F$ for the geodesic flow,
which is anything but self-explanatory, since the geodesic in~$\mathcal{F}_q$
with initial position~$q$ and initial velocity~$u_q$ is the curve given by
$\, s \mapsto \exp_Q(su_q)$, i.e., we have
\begin{eqnarray}
 &{\displaystyle
   \exp_Q(su_q) \, \Big|_{s=0}~=~q~~~,~~~
   \frac{d}{ds} \; \exp_Q(su_q) \, \Big|_{s=0}~=~u_q}& \label{eq:EXPON4}\\[2ex]
 &{\displaystyle
   \frac{D}{ds} \, \frac{d}{ds} \; \exp_Q(su_q)~=~0}& \label{eq:EXPON5}
\end{eqnarray}
Obviously, the domain $\mathrm{Dom}(\exp_Q)$ of the exponential is a tubular
neighborhood of~$Q$ in~$L\big|_Q$, and by the fundamental theorem about the
dependence of solutions of differential equations on the initial conditions
and on parameters, the map $\exp_Q$ is differentiable (i.e., smooth) and
induces the identity on~$Q$.
\begin{lem}~ \label{lem:EXPON1}
 Under the hypotheses stated above, the exponential~(\ref{eq:EXPON1}) is a
 local diffeo\-morphism onto its image, which is an open submanifold of~$P$.
\end{lem}
\proof
 As $L\big|_Q$ and $P$ have the same dimension, it suffices to prove that
 for all vectors $u_q$ in the domain of the exponential, its tangent map
 \[
  T_{u_q} \exp_Q: T_{u_q} (L\big|_Q)~\longrightarrow~T_{\exp_Q(u_q)} P
 \]
 is injective. When $u_q$ is the zero vector, this is obvious, since for all
 $\, q \in Q$, we have natural direct decompositions of the tangent spaces
 to~$L\big|_Q$ and to~$P$ at~$q$ in a `` vertical part'' and a `` horizontal
 part'',
 \[
  T_q (L\big|_Q)~=~L_q \oplus T_q Q~~~,~~~T_q P~=~L_q \oplus T_q Q
 \]
 with respect to which the tangent map
 \[
  T_q \exp_Q:~T_q(L|_Q)~\longrightarrow~T_q P
 \]
 is simply the identity. Thus let us consider the general case where
 $\, u_q \in L_q \,$ is any vector in the domain of $\exp_Q$ and
 $\, v_{u_q} \in T_{u_q} (L\big|_Q) \,$ is a tangent vector to the
 total space of the vector bundle $L\big|_Q$ over~$Q$. Suppose that
 $\, T_{u_q} \exp_Q \cdot\, v_{u_q} = 0$. Then applying the tangent
 functor to the commutative diagram
 \[
  \xymatrix{
   L\big|_Q \ar[d] \ar[r]^{\exp_Q} & P \ar[d] \\
   Q \ar[r] & P/\mathcal{F}
    }
  \] 
 where the lower horizontal arrow is the local diffeomorphism~%
 (\ref{eq:PRFOLH2}), we conclude that $v_{u_q}$ must be vertical
 and, as $\, V_{u_q} (L\big|_Q) \cong L_q$, can be identified with
 a vector $\, v_q \in L_q \,$; more explicitly, $v_{u_q} \in
 V_{u_q} (L\big|_Q) \,$ is the tangent vector
 \[
  v_{u_q}~=~\frac{d}{dt} \; (u_q + t v_q) \, \Big|_{t=0}~,
 \]
 and hence $\, T_{u_q} \exp_Q \cdot\, v_{u_q} \,$ is the tangent vector
 \[
  T_{u_q} \exp_Q \cdot\, v_{u_q}~
  =~\frac{d}{dt} \; \exp_Q(u_q + t v_q) \, \Big|_{t=0}~.
 \]
 This shows that $\, T_{u_q} \exp_Q \cdot\, v_{u_q} \,$ is the value, at
 $\, s=1$, of a Jacobi field along the geodesic $\, s \mapsto \exp_Q(su_q)$,
 defined as the variation of the following one-parameter family of geodesics:
 \[
  (s,t) \mapsto \exp_Q \bigl( s (u_q + t v_q) \bigr)~,
 \]
 where~$t$ is the family parameter and $s$ is the geodesic parameter
 (for fixed~$t$).
 Explicitly, the value of this Jacobi field at the point $\exp_Q(su_q)$ is
 \[
  \frac{d}{dt} \; \exp_Q \bigl( s (u_q + t v_q) \bigr) \, \Big|_{t=0}~,
 \]
 showing that if $\, T_{u_q} \exp_Q \cdot\, v_{u_q} = 0$, it must vanish
 at $\, s=0 \,$ and at $\, s=1$, i.e., $q$ and $\, p = \exp_Q(u_q) \,$
 would be conjugate points along the geodesic $\, s \mapsto \exp_Q(su_q)$.
% as shown in Figure~\ref{fig:JAC}.
% \begin{figure}[ht]
% \begin{center}
%  \psfrag{0}{$0$}
%  \psfrag{u}{$u_q$}
%  \psfrag{v}{$v_q$}
%  \psfrag{L}{$L_q$}
%  \psfrag{EX}{$\exp_Q$}
%  \psfrag{F}{$\mathcal{F}_q$}
%  \psfrag{q}{$q$}
%  \psfrag{p}{$p$}
%  \psfrag{X}{$T_{u_q} \exp_Q \cdot\, v_{u_q}$}
%  \includegraphics[width=10cm,height=8cm]{figJACOB.eps}
%  \caption{Zeros of the derivative of the exponential and Jacobi fields}
%  \label{fig:JAC}
% \end{center}
% \end{figure}
%
% \noindent
 But the condition that the connection $\nabla$ has vanishing curvature
 and torsion excludes the existence of conjugate points along any geodesic,
 because the differential equation for a Jacobi field, written in components
 for an autoparallel frame along that geodesic, reduces to an equation of
 the form $\, d^{\,2} X^i/ds^2 = 0$, whose solutions have exactly one
 zero~-- no more, no less.
\qed

\noindent
In particular, it follows that the exponential~(\ref{eq:EXPON1}) provides
a diffeomorphism
\begin{equation} \label{eq:EXPON6}
 \exp_Q: U_0~\longrightarrow U
\end{equation}
of a convex neighborhood $U_0$ of $Q$ in the vector bundle $L\big|_Q$
with a tubular neighborhood $U$ of~$Q$ in~$P$.%
\footnote{A neighborhood of the zero section of a vector bundle is called
convex if its intersection with each fiber is a convex neighborhood of the
origin in that fiber.}
By construction, this diffeomorphism is affine.

In the case of geodesic completeness, we can prove an even stronger claim:
\begin{lem}~ \label{lem:EXPON2}
 Under the hypotheses stated above, and if $P$ is geodesically complete
 with respect to~$\nabla$, the exponential~(\ref{eq:EXPON1}) defines a
 covering
 \begin{equation} \label{eq:EXPON7}
  \exp_Q:~L\big|_Q~\longrightarrow~P
 \end{equation}
 which, for every point $q$ of~$Q$, induces a universal covering
 \begin{equation} \label{eq:EXPON8}
  \exp_q:~L_q~\longrightarrow~\mathcal{F}_q
 \end{equation}
 of the leaf $\mathcal{F}_q$ by the fiber $L_q$. In particular, if all
 leaves $\mathcal{F}_q$ are simply connected, the exponential provides
 a global affine diffeomorphism between $L\big|_Q$ and~$P$.
\end{lem}
\proof
 Under the hypothesis of geodesic completeness, the domain of the exponential
 $\exp_Q$ is the entire vector bundle $L\big|_Q$. Moreover, the hypothesis
 that the connection $\nabla$ should have vanishing curvature and torsion
 implies that for every point $q$ of~$Q$, the restriction $\exp_q$ of
 $\exp_Q$ to the fiber $L_q$ is an affine map from the vector space~$L_q$,
 equipped with the trivial linear connection, to the leaf~$\mathcal{F}_q$,
 equipped with the linear connection $\, \nabla_q = \nabla\big|_{\mathcal{F}_q}$.
 (For a much more general statement, see, for example, \cite[Chapter~6,
 Theorem~7.1, p.~257]{KN1}.) So, the lemma follows from a general theorem,
 stated in more detail and proved in Appendix~B, according to which every
 affine map from a connected, simply connected and geodesically complete
 affine manifold $M$ to a connected affine manifold $M'$, if it is a local
 diffeomorphism, it is even a covering; in particular, it is automatically
 surjective (and $M'$ is automatically geodesically complete).
\qed

Another result which we shall need in what follows concerns differential forms
on foliated manifolds:
\begin{prp}~\label{prp:PULLB}
 Let\/ $P$ be a manifold equipped with a simple foliation\/~$\mathcal{F}$
 with involutive tangent distribution\/~$L$, and let\/ $\alpha$ be a\/
 $k$-form on\/~$P$. Then\/ $\alpha$ is the pull-back of a\/~$k$-form\/
 $\alpha_Q^{}$ on the quotient manifold\/~$\, Q = P/\mathcal{F} \,$ by the
 projection\/ $\, \pi: P \longrightarrow Q$, $\alpha = \pi^* \alpha_Q^{}$,
 if and only if, for all vector fields\/~$X$ along\/~$\mathcal{F}$,
 $X \in \Gamma(L)$, we have $\, i_X^{} \alpha = 0$ (horizontality
 condition) and $\, \mathbb{L}_X^{} \alpha = 0$ (condition of
 constancy along the leaves).
 \end{prp}

\proof
 First, observe that for $\alpha$ to be the pull-back of a $k$-form
 $\alpha_Q^{}$ on the quotient manifold~$Q$, we must have
 \begin{equation} \label{eq:PULLB}
  \begin{array}{c}
   \alpha(p) \bigl( u_1,\ldots,u_k \bigr)~
   =~\alpha_Q(\pi(p))
     \bigl( T_p \pi \cdot u_1 , \ldots , T_p \pi \cdot u_k \bigr) \\[1mm]
   \mbox{for $\, p \in P$, $u_1,\ldots,u_k \in T_p P$}
  \end{array}~.
 \end{equation}
 Thus if any of the vectors $u_1,\ldots,u_k$ belongs to~$L_p$, this
 expression vanishes, so for any $\, X \in \Gamma(L)$, we must
 have $\, i_X^{} \alpha = 0 \,$ and therefore also
 \[
  \mathbb{L}_X^{} \alpha~
  =~\bigl( d i_X^{} + i_X^{} d \bigr) \, \pi^* \alpha_Q^{}~
  =~i_X^{} \pi^* (d \alpha_Q^{})~=~0~.
 \]
 Conversely, it is clear that if we use eqn~(\ref{eq:PULLB}) to define
 $\alpha_Q^{}$ in terms of~$\alpha$, we must ensure (i)~that for fixed
 $\, p \in P$, the expression on the lhs of this equation does not
 depend on the representatives $\, u_i \in T_p P \,$ of the vectors
 $\,  T_p \pi \cdot u_i \in T_{\pi(p)} Q$, which is guaranteed by the
 condition of horizontality ($i_X^{} \alpha = 0 \,$  for $\, X \in
 \Gamma(L)$), and (ii)~that the expression on the lhs of this equation 
 does not depend on the representative $\, p \in P \,$ of the point
 $\, \pi(p) \in Q \,$: this can be derived from the condition of
 constancy along the leaves ($\mathbb{L}_X^{} \alpha = 0 \,$ for
 $\, X \in \Gamma(L)$), as follows: Let $p$ and $p'$ be two points
 of~$P$ such that $\, \pi(p) = \pi(p') \,$: this means that they
 belong to the same leaf~$F$, and therefore there is a curve $\gamma$
 entirely contained in the leaf~$F$ with $\, \gamma(0) = p \,$ and
 $\, \gamma(1) = p' \,$; in particular, we have $\, \dot{\gamma}(s)
 \in L_{\gamma(s)}$ \linebreak for $\, 0 \leqslant s \leqslant 1$,
 and we can further assume that $\, \dot{\gamma}(s) > 0 \,$ for
 $\, 0 \leqslant s \leqslant 1$. Using a partition $(s_\alpha)_%
 {\alpha=1,\ldots,r}$ of the interval $[0,1]$ ($0 = s_0 < s_1 <
 \ldots < s_r < s_{r+1} = 1$) and a finite family $(U_\alpha)_%
 {\alpha=0,\ldots,r}$ of chart domains $U_\alpha$ for~$P$ where
 $\pi$ is represented by a constant projection onto some subspace,
 together with some smooth cutoff function of compact support on~$P$
 that is $1$ on an open neighborhood of the image of the curve $\gamma$,
 it becomes evident that we can find a vector field $X$ on~$P$ along~%
 $\mathcal{F}$, $X \in \Gamma(L)$, which extends $\dot{\gamma}$, i.e.,
 such that $\, \dot{\gamma}(s) = X(\gamma(s)) \,$ for $\, 0 \leqslant
 s \leqslant 1$. But this means that $\gamma$ is an integral curve
 of~$X$ and, more than this, that close to $\, p = \gamma(0)$, the
 flow $F_X^{}$ of~$X$ is defined at least up to $\, s=1$, so that
 there exist open neighborhoods $U_0$ of $\, p = \gamma(0) \,$ and
 $U_1$ of $\, p' = \gamma(1) \,$ such that the flow for time~$1$
 establishes a diffeomorphism $\, F_X^{}(1,.): U_0 \longrightarrow U_1$
 which preserves the leaves of~$\mathcal{F}$, since $X$ is tangent
 to~$\mathcal{F}$, i.e., we have $\, \pi\big|_{U_1} = F_X^{}(1,.)
 \smcirc \pi\big|_{U_0}$. Now, $\mathbb{L}_X^{} \alpha = 0 \,$
 implies $\, F_X(1,.)^* (\alpha\big|_{U_1}) = \alpha\big|_{U_0}$,
 so using $\, p' = F_X^{}(1,p) \,$ and setting $\, u'_i = T_p^{}
 F_X^{}(1,.) \cdot u_i$ ($1 \leqslant i \leqslant k$), we obtain
 $\, T_{p'} \pi \cdot u'_i = T_p \pi \cdot u_i$ ($1 \leqslant i
 \leqslant k$) and
 \begin{eqnarray*}
  \alpha_{p'}^{} \bigl( u'_1,\ldots,u'_k \bigr) \!\!
  &=& \alpha_{F_X(1,p)}^{} \bigl( T_p^{} F_X^{}(1,.) \cdot u_1^{} , \ldots ,
                              T_p^{} F_X^{}(1,.) \cdot u_k^{} \bigr) \\
  &=& \bigl( F_X(1,.)^* (\alpha\big|_{U_1}) \bigr)_p^{}
      \bigl( u_1^{} , \ldots , u_k^{} \bigr) \\
  &=& \alpha_p^{} \bigl( u_1^{} , \ldots , u_k^{} \bigr)~.
 \end{eqnarray*}
\qed

\section{Foliated symplectic manifolds}

In this section, we consider the geometry of a foliated symplectic
manifold, or more precisely, of a symplectic manifold~$P$, with
symplectic form~$\,\omega$, that comes equipped with a simple
lagrangian foliation~$\mathcal{F}$.
Of course, lagrangian foliations may exist or not, and they can
be simple or not: a classical example of a symplectic manifold
which does not admit any lagrangian foliation is the sphere $S^2$
(``no-hair theorem''), while a classical example of a lagrangian
foliation which is regular but not simple is the irrational flow
on the torus~$\mathbb{T}^2$ (in both cases, the symplectic form
is the standard volume form).
But in the case of a simple lagrangian foliation, the quotient
space $\, Q = P/\mathcal{F} \,$ admits a unique manifold structure
such that the canonical projection $\pi$ from~$P$ to $Q$ is a
surjective submersion.
Note that with this convention, $Q$ is a quotient manifold of~$P$,
but nothing guarantees ``a~priori'' that it can be realized as a
submanifold of~$P$, so the existence of an embedding of~$Q$
into~$P$ as a closed submanifold is an additional condition
that must be imposed separately or deduced from other
additional assumptions.%
\footnote{In general, there may be topological obstructions to the
existence of an embedding of~$Q$ into~$P$. Such obstructions are of
global nature, since the submersion theorem (or local slice theorem)
states that locally, there is always such an embedding. As an example
of a set of additional conditions that guarantees its global existence,
we mention the hypotheses in the third item of Theorem~\ref{thm:PULLB}
below~-- namely that the leaves are geodesically complete with respect
to the Bott connection and simply connected, as these ensure that $P$
is an affine bundle over~$Q$, and affine bundles always admit global
sections.}
In any case, the hypothesis that the foliation $\mathcal{F}$ is
lagrangian provides a canonical partial linear connection in~$L$
along~$L$ with vanishing curvature and torsion, namely the Bott
connection $\nabla^B$ introduced earlier: it implies that the
leaves of~$\mathcal{F}$ are flat affine manifolds and is the
crucial ingredient in the proof of the following statement:
\begin{thm}~ \label{thm:PULLB}
 Let\/ $P$ be a symplectic manifold, with symplectic form~$\,\omega$,
 equipped with an involutive lagrangian distribution\/~$L$. Suppose
 that the corresponding foliation\/~$\mathcal{F}$ is simple, writing
 its leaves as the level sets of a surjective submersion\/ $\, \pi_P:
 P  \longrightarrow Q$, and suppose finally that the quotient manifold\/
 $\, Q = P/\mathcal{F} \,$ can be realized as a closed embedded sub%
 manifold of\/~$P$. Under these circumstances, consider the musical
 isomorphism \linebreak $\, \omega^\sharp: L\big|_Q^\perp \longrightarrow
 L\big|_Q$ (see eqn~(\ref{eq:MUSISO4})), together with the isomorphism\/
 $\, L\bigl|_Q^\perp ~\cong\, T^* Q \,$ (which arises from the direct
 decomposition~(\ref{eq:SVCOMP1})), and combined with the exponential\/~%
 $\exp_Q$ as defined in Sect.~7. Then we have the following:
 \begin{itemize}
  \item The composition of these isomorphisms provides a diffeomorphism\/
        $\, \phi: V \longrightarrow U$ \linebreak of a tubular neighborhood\/
        $U$ of\/~$Q$ in\/~$P$ with a convex neighborhood\/~$V$ of the zero
        section of the cotangent bundle\/ $T^* Q$ of\/~$Q$.
  \item If the leaves of\/~$\mathcal{F}$ are geodesically complete with
        respect to the Bott connection, the composition of these isomorphisms
        provides a covering of\/~$P$ by the cotangent bundle of\/~$Q$,
        $\phi: T^* Q \longrightarrow P \,$.
  \item If the leaves of\/~$\mathcal{F}$ are geodesically complete with
        respect to the Bott connection and simply connected, the composition
        of these isomorphisms provides a diffeomorphism of\/~$P$ with the
        cotangent bundle of\/~$Q$, $\phi: T^* Q \longrightarrow P \,$.
 \end{itemize}
 Furthermore, $\phi$ preserves fibers, mapping\/~$T_q^* Q$ onto\/~%
 $\mathcal{F}_q^{}$ (or, in the first case, $V \cap  T_q^* Q$ \linebreak
 onto\/ $\, U \cap \mathcal{F}_q^{}$), and defining
 \begin{equation} \label{eq:1FORM}
  \theta~= \; - \, i_\Sigma^{} \, \phi^*\omega~,
 \end{equation}
 we have that $\, \phi^*\omega + d\theta \,$ is the pull-back of a
 closed\/~$2$-form $\,\omega_Q^{}$ on\/~$Q$ by the projection\/~$\tau$
 of\/~$T^*Q$ to\/~$Q$:
 \begin{equation}\label{eq:PULLB1}
  \phi^*\omega + d \theta~=~\tau^* \omega_Q^{}~.
 \end{equation}
 Finally, the  cohomology class $\, [\,\omega_Q^{}] \in H^2(Q) \,$ of~%
 $\,\omega_Q^{}$ does not depend on the embedding employed and thus is
 an invariant of the foliation\/~$\mathcal{F}$.
\end{thm}
\begin{rmk}~
 The last statement ensures that $\phi$ is ``almost'' a symplectomorphism:
 $\phi^* \omega$ differs from the standard symplectic form of the cotangent 
 bundle only by the pull-back of a closed $2$-form on the base. If $Q$ is
 a lagrangian submanifold of~$P$, then $\ \omega_Q^{} = 0$ \linebreak
 and $\phi$ will be a symplectomorphism. In this special case, the first
 statement of the above theorem, which is of local nature (with respect to
 the structure of~$P$ along the leaves of the foliation~$\mathcal{F}$),
 is known as Weinstein's symplectic tubular neighborhood theorem,
 established in~\cite{We}. The third statement has first been proved
 in~\cite{Th}. Here, besides establishing also the second statement,
 we give a more direct proof for all three of them, avoiding the use
 of additional and artificial ingredients (such as the auxiliary
 riemannian metric employed in~\cite{Th}): this will also allow us
 to formulate and prove an extension of this theorem to the case of
 polysymplectic and multisymplectic geometry, treated in the next
 section.
\end{rmk}
\proof
 In view of Lemmas~\ref{lem:EXPON1} and~\ref{lem:EXPON2}, we just
 need to prove the final part, contained in eqns~(\ref{eq:1FORM})
 and~(\ref{eq:PULLB1}).
 To simplify the presentation, we consider only the first and third
 statement, where $\phi$ is a diffeomorphism and hence can be used
 to identify $V$ with~$U$ and $T^* Q$ with~$P$, respectively.
 (The second statement, where $\phi$ is just a local diffeomorphism,
 can be treated similarly, taking into account that in this case,
 the Euler vector field $\Sigma$ may fail to be globally defined
 on~$P$, but it can be replaced by a family of Euler vector fields
 locally defined on~$P$, which leads to a family of local formulas
 of the same type as eqns~(\ref{eq:1FORM}) and~(\ref{eq:PULLB1}).)
 Therefore, we suppress the reference to the pull-back by~$\phi$.
 The argument will be based on Proposition~\ref{prp:PULLB}, according
 to which it is sufficient to show that for every vertical vector
 field~$X$, we have $\, i_X^{}(\omega + d\theta) = 0$ (horizontality
 condition) and $\, \mathbb{L}_X^{}(\omega + d\theta) = 0$ (condition
 of constancy along the leaves). Since $\,\omega$ is closed, the
 second of these conditions follows directly from the first:
 \[
  i_X^{}(\omega+d\theta)~=~0~~~\Longrightarrow~~~
  \mathbb{L}_X^{}(\omega + d\theta)~
  =~\bigl( d i_X^{} + i_X^{} d \bigr)(\omega + d\theta)~
  =~d \bigl( i_X^{}(\omega+d\theta) \bigr)~=~0~.
 \]
 To prove the first, we must show that for every vector
 field~$X$ along~$\mathcal{F}$ and every vector field~$Y$,
 we have $\; (\omega + d\theta)(X,Y) = 0$, and due to
 Lemma~\ref{lem:CVPROJ2}, we may do so assuming, without
 loss of generality, that $Y$ is projectable. Now using the
 definitions of the Bott connection and of the Euler vector
 field, we have
 \[
  \omega(X,Y)~=~\omega(\nabla^{B}_{X}\Sigma, Y)~
  =~X \!\cdot \omega(\Sigma,Y) \, - \, \omega(\Sigma,[X,Y])~.
 \]
 Since $X$ is vertical and $Y$ is projectable, $[X,Y]$ is also
 vertical (see eqn~(\ref{eq:CVPROJ2})), and since $L$ is lagrangian,
 the second term vanishes, so we get
 \[
  \omega(X,Y)~=~X \!\cdot \omega(\Sigma,Y)~.
 \]
 Using the definition of~$\theta$, eqn~(\ref{eq:1FORM}), together
 with the fact that this implies that $\theta$ vanishes on vertical
 vector fields, we have
 \begin{eqnarray*}
  (\omega+d\theta)(X,Y) \!\!
  &=&\!\! \omega(X,Y) + X \!\cdot \theta(Y)
          - Y \!\cdot \theta(X) - \theta([X,Y]) \\
  &=&\!\! X \!\cdot \omega(\Sigma,Y) + X \!\cdot \theta(Y)
          - Y \!\cdot \theta(X) - \theta([X,Y]) \\
  &=&\!\! \mbox{} - X \!\cdot \theta(Y) + X \!\cdot \theta(Y)
          - Y \!\cdot \theta(X) - \theta([X,Y]) \\
  &=&\!\! 0~.
 \end{eqnarray*}
 Finally, we must address the issue of uniqueness, or rather the amount of
 non-uniqueness, of the decomposition~(\ref{eq:PULLB1}), generated by the
 fact that there are different Euler vector fields, corresponding to different
 choices of the embedding of the quotient manifold $P/\mathcal{F}$ into~$P$.
 Thus let $\Sigma_1$ and $\Sigma_2$ be two Euler vector fields, and define
 $\, \theta_1^{} = - i_{\Sigma_1} \omega \,$ and $\, \theta_2^{} =
 - i_{\Sigma_2} \omega$. Then for every vertical vector field~$X$,
 \[
  i_X^{} (\theta_1-\theta_2)~=~\omega(X,\Sigma_2-\Sigma_1)~=~0~,
 \]
 since $L$ is lagrangian, while for every projectable vector field~$Y$,
 \begin{eqnarray*}
  \mathbb{L}_X^{} (\theta_1-\theta_2) \, (Y) \!\!
  &=&\!\! X \cdot\, ((\theta_1-\theta_2)(Y)) \, - \,
          (\theta_1-\theta_2)([X,Y]) \\
  &=&\!\! X \cdot\, \omega(\Sigma_2-\Sigma_1,Y) \, - \,
          \omega(\Sigma_2-\Sigma_1,[X,Y]) \\
  &=&\!\! \omega(\nabla_X^B(\Sigma_2-\Sigma_1),Y) \\
  &=&\!\! 0~,
 \end{eqnarray*}
 where we have used the definition of the Bott connection and the fact
 that $\, \Sigma_2 - \Sigma_1 \,$ is covariantly constant. According
 to Proposition~\ref{prp:PULLB}, it follows that there is a $1$-form
 $\theta_Q$ on~$Q$ such that $\, \theta_1 - \theta_2 =~\pi^* \theta_Q$,
 implying
 \[
  \omega + d \theta_1^{}~=~\pi^* \omega_Q^{(1)}~~,~~
  \omega + d \theta_2^{}~=~\pi^* \omega_Q^{(2)}~,
 \]
 with
 \[
  \theta_1^{} - \theta_2^{}~=~\pi^* \theta_Q^{}~~,~~
  \omega_Q^{(1)} - \omega_Q^{(2)}~=~d \theta_Q^{}~.
  \vspace{1ex}
 \]
 In particular, the cohomology class of~$\omega_Q^{}$ does not depend
 on the choice of embedding.
\qed

\section[Structure of polysymplectic and multisymplectic fiber bundles]
{Structure of polysymplectic and \\ multisymplectic fiber bundles}

In analogy with the symplectic case, we can now formulate our main theorem
about the structure of polysymplectic and multisymplectic fiber bundles.
The additional ingredient, as compared to the symplectic case, comes from
the fact that the underlying manifold~$P$ is now the total space of a
fiber bundle over some other manifold\/~$M$,%
\footnote{In applications to physics, the base manifold~$M$ is space-time.
In classical mechanics, this reduces to a copy of the real line (time axis)
which is usually suppressed, but it reappears immediately when one considers
non-autonomous systems, passing from symplectic manifolds to contact manifolds
and then, using Cartan's trick of adding yet another copy of the real line
(energy axis), back to symplectic manifolds.}
with bundle projection denoted by $\, \pi: P \longrightarrow M$, and
that the distribution~$L$ is vertical with respect to this projection.
Roughly speaking, this implies that the submanifold of~$P$ representing
the quotient space $P/\mathcal{F}$, which is now denoted by~$E$, should be
the total space of a fiber bundle on~$M$, whose projection will be denoted
by $\, \pi_E : E \longrightarrow M$, as in the examples in Sects~4 and~5.
Thus, the condition that $E$ is a submanifold of~$P$ complementary to~%
$\mathcal{F}$ and, at the same time, to the fibers of the projection~$\pi$,
leads us to replace eqn~(\ref{eq:SVCOMP1}) by the condition that for every
point $e$ of~$E$, we have
\begin{equation} \label{eq:SVCOMP2}
 T_e P~=~T_e E \,\oplus\, L_e \qquad \mbox{and} \qquad
 V_e P~=~V_e E \,\oplus\, L_e~,
\end{equation}
where $V_e P$ denotes the vertical space with respect to the projection $\pi$
and $V_e E$ denotes the vertical space with respect to the projection $\pi_E$.

In the case of polysymplectic fiber bundles, we have
\begin{thm}~ \label{thm:POLPULLB} 
 Let\/ $P$ be a polysymplectic fiber bundle over a manifold\/~$M$,
 with projection \linebreak $\pi: P \longrightarrow M$, polysymplectic
 form\/~$\,\hat\omega$ and involutive polylagrangian distribution\/~$L$.
 \linebreak
 Suppose that the corresponding foliation\/~$\mathcal{F}$ is simple,
 writing its leaves as the level sets of a surjective submersion\/
 $\, \pi_P: P \longrightarrow E$, that $\pi$ induces a surjective
 submersion\/ $\,\pi_E: E \longrightarrow M \,$ so that\/ $\, \pi
 = \pi_E \smcirc \pi_P$, and finally that the quotient manifold\/
 $\, E = P/\mathcal{F}$ \linebreak (a)~can be realized as a closed
 embedded submanifold of\/~$P$ and (b)~is the total space of a fiber
 bundle over\/~$M$ with respect to the projection\/~$\pi_E$.
 Under these circumstances, consider the musical isomorphism
 $\; \hat{\omega}^\sharp: \bwedge{k} L\big|_E^\perp \otimes\, \pi_E^*(\hat{T})
 \longrightarrow L\big|_E$ (see eqn~(\ref{eq:MUSISO5})), together with
 the isomorphism $\, L\bigl|_E^\perp ~\cong \, V^* E \,$ (which arises
 from the direct decomposition~(\ref{eq:SVCOMP2})), and combined with
 the exponential\/~$\exp_E$ as defined in Sect.~7.
 Then we have the following:
 \begin{itemize}
  \item The composition of these isomorphisms provides a diffeomorphism\/
        $\, \phi: V \longrightarrow U$ \linebreak of a tubular neighborhood\/
        $U$ of\/~$E$ in\/~$P$ with a convex neighborhood\/~$V$ of the zero
        section of the model vector bundle\/ $\, \bwedge{k} V^* E \otimes
        \pi_E^*(\hat{T}) \,$ of Example~\ref{exe:MPFPL}.
  \item If the leaves of\/~$\mathcal{F}$ are geodesically complete with
        respect to the Bott connection, the composition of these isomorphisms
        provides a covering of\/~$P$ by the model vector bundle, $\, \phi:
        \bwedge{k} V^* E \otimes \pi_E^*(\hat{T}) \longrightarrow P \,$.
  \item If the leaves of\/~$\mathcal{F}$ are geodesically complete with
        respect to the Bott connection and simply connected, the composition
        of these isomorphisms provides a diffeomorphism of\/~$P$ with the
        model vector bundle, $\, \phi: \bwedge{k} V^* E \otimes
        \pi_E^*(\hat{T}) \longrightarrow P \,$.
 \end{itemize}
 Furthermore, $\phi$ preserves fibers, mapping\/~$\, \bwedge{k} V_e^* E
 \otimes \hat{T}_{\pi_E(e)^{}} \,$ onto\/~$\mathcal{F}_e^{}$ (or, in the
 first case, $\, V \cap\, \bwedge{k} V_e^* E \otimes \hat{T}_{\pi_E(e)} \,$
 onto\/ $\, U \cap\, \mathcal{F}_e^{}$), and defining
 \begin{equation} \label{eq:kFORM1}
  \hat{\theta}~= \; - \, i_\Sigma^{} \, \phi^* \hat{\omega}~,
 \end{equation}
 we have that $\, \phi^* \hat{\omega} + d_V \hat{\theta} \,$ is the pull-back
 of a vertically closed\/~$k$-form $\,\hat{\omega}_E^{}$ on\/~$E$ by the
 projection\/~$\pi^k$ of\/~$\, \bwedge{k} V^* E \otimes \pi^*(\hat{T}) \,$
 to\/~$E$:
 \begin{equation}\label{eq:POLPULLB1}
  \phi^* \hat{\omega} + d_V \hat{\theta}~=~(\pi^k)^* \hat{\omega}_E^{}~.
 \end{equation}
 Finally, the cohomology class $\, [\,\hat{\omega}_E^{}] \in H^k(E) \,$
 of\/~$\,\hat\omega_E^{}$ does not depend on the embedding employed and
 thus is an invariant of the foliation\/~$\mathcal{F}$.
\end{thm}
\proof
 The proof is completely analogous to the proof of Theorem~\ref{thm:PULLB}
 for foliated symplectic manifolds, and the calculations to verify the formula 
 (\ref{eq:kFORM1}) are carried out with a vector field~$X$ along~$\mathcal{F}$
 and $k$ projectable vector fields $Y_1,\ldots,Y_k$, all of them vertical with
 respect to the projection $\pi$ to~$M$.
\qed 

Turning to the case of multisymplectic fiber bundles, we have
\begin{thm}~ \label{thm:MULPULLB}
 Let\/ $P$ be a multisymplectic fiber bundle over a manifold\/~$M$,
 with projection \linebreak $\pi: P \longrightarrow M$, multisymplectic 
 form\/~$\,\omega$ and involutive multilagrangian distribution\/~$L$.
 \linebreak
 Suppose that the corresponding foliation\/~$\mathcal{F}$ is simple,
 writing its leaves as the level sets of a surjective submersion
 $\, \pi_P: P \longrightarrow E$, that\/ $\pi$ induces a surjective
 submersion\/ $\, \pi_E: E \longrightarrow M \,$ so that $\, \pi
 = \pi_E \smcirc \pi_P$, and finally that the quotient manifold\/
 $\, E = P/\mathcal{F}$ \linebreak (a)~can be realized as a closed
 embedded submanifold of\/~$P$ and (b)~is the total space of a fiber
 bundle over\/~$M$ with respect to the projection\/~$\pi_E$.
 Under these circumstances, consider the musical isomorphism
 $\; \omega^\sharp: \dbwedge{\,r-1}{\;~k} L\big|_E^\perp
 \longrightarrow L\big|_E$ (see eqn~(\ref{eq:MUSISO6})),
 together with the isomorphism $\, L\bigl|_E^\perp ~\cong\, T^* E \,$
 (which arises from the direct decomposition~(\ref{eq:SVCOMP2})), and
 combined with the exponential\/~$\exp_E$ as defined in Sect.~7.
 Then we have the following:
 \begin{itemize}
  \item The composition of these isomorphisms provides a diffeomorphism\/
        $\, \phi: V \longrightarrow U$ \linebreak of a tubular neighborhood\/
        $U$ of\/~$E$ in\/~$P$ with a convex neighborhood\/~$V$ of the zero
        section of the model vector bundle\/ $\, \dbwedge{\,r-1}{\;~k}
        T^* E \,$ of Example~\ref{exe:MPFML}.
  \item If the leaves of\/~$\mathcal{F}$ are geodesically complete with
        respect to the Bott connection, the composition of these isomorphisms
        provides a covering of\/~$P$ by the model vector bundle, $\, \phi:
        \dbwedge{\,r-1}{\;~k} T^* E \longrightarrow P \,$.
  \item If the leaves of\/~$\mathcal{F}$ are geodesically complete with
        respect to the Bott connection and simply connected, the composition
        of these isomorphisms provides a diffeomorphism of\/~$P$ with the
        model vector bundle, $\, \phi: \dbwedge{\,r-1}{\;~k} T^* E
        \longrightarrow P \,$.
 \end{itemize}
 Furthermore, $\phi$ preserves fibers, mapping\/~$\, \dbwedge{\,r-1}{\;~k}
 T_e^* E \,$ onto\/~$\mathcal{F}_e^{}$ (or, in the first case, \linebreak
 $\, V \cap \dbwedge{\,r-1}{\;~k} T_e^* E \,$ onto\/
 $\, U \cap \mathcal{F}_e^{}$), and defining
 \begin{equation} \label{eq:kFORM2}
  \theta~= \; - \, i_\Sigma^{} \, \phi^* \omega~,
 \end{equation}
 we have that $\, \phi^* \omega + d\theta \,$ is the pull-back of a
 closed\/~$k$-form $\,\omega_E^{}$ on\/~$E$ by the projection\/~%
 $\pi_{r-1}^{~k}$ of\/~$\, \dbwedge{\,r-1}{\;~k} T^* E \,$ to\/~$E$:
 \begin{equation}\label{eq:MULPULLB1}
  \phi^*\omega + d\theta~=~(\pi_{r-1}^{~k})^* \omega_E^{}~.
 \end{equation}
 Finally, the cohomology class $\, [\,\omega_E^{}] \in H^k(E) \,$ of\/~%
 $\,\omega_E^{}$ does not depend on the embedding employed and thus is an
 invariant of the foliation\/~$\mathcal{F}$.
\end{thm}
\proof 
 Analogous to that of Theorem~\ref{thm:POLPULLB}, eliminating only the 
 condition of verticality of the projectable vector fields relative to
 the projection over~$M$.
\qed

\section{Conclusions}

The main new results reported in this paper are the theorems on existence
of torsion-free polysymplectic and multisymplectic connections and their
complete classification (Theorems~\ref{thm:PLCONN1} and~\ref{thm:MLCONN1}),
together with the structure theorems on polysymplectic and multisymplectic
fiber bundles (Theorems~\ref{thm:POLPULLB} and~\ref{thm:MULPULLB}) which
show that, under certain mild additional assumptions, these are exhausted
by the well-known standard examples of bundles of forms (Examples~%
\ref{exe:MPFPL} and~\ref{exe:MPFML}).
All these generalize corresponding theorems of symplectic geometry which
we have decided to include not only for the sake of completeness (given
that most of them do not appear to have been stated explicitly in the
existing literature, at least not in their full generality), but also
because our proofs use different techniques.
For example, we could not find an explicit statement of the theorem on
the existence and classification of torsion-free symplectic connections
that preserve a single lagrangian foliation, included here as Theorem~%
\ref{thm:LAGCON}.
(What one can find easily are classification theorems for torsion-free
symplectic connections which either are subject to no further constraints
or else are required to preserve two transversal lagrangian foliations:
as is well known, the latter case leads to a unique answer, namely the
bilagrangian connection first constructed by He\ss.
This situation is somewhat surprising since after all, the case of a
single \mbox{lagrangian} foliation \emph{is} very important: it is
the situation one encounters
when dealing with cotangent bundles! Indeed, a cotangent bundle is a
symplectic manifold carrying a distinguished lagrangian foliation but
no natural candidate for a second one that would be transversal to it:
all one finds are single lagrangian submanifolds transversal to it,
namely the zero section or, more generally, the graph of any closed
$1$-form on the base manifold.)
Similarly, the global versions of Weinstein's tubular neighborhood
theorem do not seem to have been formulated in their full generality,
and the existing proofs use rather artificial additional ingredients
which, as we show, are really unnecessary.

\pagebreak

Regarding the extension from symplectic to polysymplectic and
multisymplectic geo\-metry, one of the central concepts is the
Bott connection: it is a partial linear connection in and along
the corresponding polylagrangian or multilagrangian distribution~$L$
and is a natural geometric object at least when $L$ is uniquely
determined and involutive, which is the generic case~\cite{FG}.
Since this connection is both torsion-free and flat, it implies
that, just as in symplectic geometry, the leaves of the
corresponding foliation are flat affine manifolds~-- a fact
that imposes severe restrictions on the underlying geometry.
\linebreak
\emph{The upshot is that polysymplectic and multisymplectic
geometry is analogous not to the geometry of general symplectic
manifolds but rather to that of foliated symplectic manifolds,
and that is why there is no generic polysymplectic or multi%
symplectic analogue of the coadjoint orbit construction, since
typically such orbits do not admit lagrangian foliations.}

Such observations, when applied to classical field theory,
support a general picture concerning the role of position
variables and momentum variables in physics.

In the usual hamiltonian formulation of classical mechanics,
these variables are essentially treated on an equal footing:
they can be thought of as ingredients of local coordinate
systems in a symplectic manifold, called phase space, and
transformations between such local coordinate systems, called
canonical transformations, are symmetries of the theory~--
a point of view that has been triumphant in the mathematical
treatment of completely integrable systems, whose solution is
achieved through a judiciously chosen canonical transformation
to so-called action-angle variables.
As a result, many have been led to believe that there is a
general ``democracy'' between position and momentum variables.

However, it is well known that this ``democracy'' is lost upon
quantization: in contrast to what happens in classical mechanics,
canonical transformations mixing position and momentum variables
are no longer symmetries of quantum mechanics, since they cannot
be implemented by unitary operators in the Hilbert space of states.

What is much less known is that this loss of symmetry is by no
means a specific feature of going to the quantum world, simply
because the same thing happens in (relativistic) field theory:
here too, this ``democracy'' just disappears!

The central reason seems to be that, already at the classical
level, relativity is built on fundamental new principles of
physics that require a clear-cut distinction between the two
types of variables.
Perhaps the most important of them all is space-time locality,
which postulates that events localized in space-like separated
regions of space-time cannot exert any direct influence on each
other: obviously, this principle refers to space-like separation
in space-time and not in momentum space!
Therefore, it is not a defect but rather a virtue of polysymplectic
and multisymplectic geometry, whose proposal is to provide the correct
mathematical framework for the hamiltonian formulation of (relativistic)
classical field theory, that they incorporate, from the very beginning,
a clear geometrical distinction between position and momentum variables,
in terms of a given distribution describing the ``collection of all
momentum directions''.
This characterization is as it should be: coordinate and frame
independent, as well as intrinsically defined and unique; its
 mere existence being in sharp contrast to the situation in
symplectic geometry, where specifying a lagrangian distribution
is a matter of choice.
Thus in (relativistic) field theory, the lack of ``democracy''
in the sense described before is \emph{not} a quantum effect,
but rather the result of physical principles which already
prevail at the classical level.

\begin{appendix}

\section*{Appendix A: Auxiliary formulas}

In the course of this paper, we have repeatedly made use of the following
two elementary facts.
\begin{lem}~ \label{lem:CONN1}
 Let\/ $M$ be a manifold and let\/ $L$ be a distribution on\/~$M$.
 Then if there exists a torsion-free linear connection\/ $\nabla$
 on\/~$M$ preserving\/~$L$, or more generally, if there exist an
 involutive distribution\/~$V$ on~$M$ containing\/~$L$ and a
 torsion-free partial linear connection\/ $\nabla$ in\/~$V$
 along\/~$V$ preserving\/~$L$, $L$ must be involutive.
\end{lem}
\proof
 This follows simply by looking at the definition of the torsion tensor
 of~$\nabla$,
 \[
  T(X,Y)~=~\nabla_{\!X}^{} Y \, - \, \nabla_{\!Y}^{} X \, - \, [X,Y]
 \]
 which implies that if $\nabla$ is torsion-free and preserves~$L$, then
 when $X$ and~$Y$ are along~$L$, so must be~$[X,Y]$.
\qed
\begin{lem}~ \label{lem:CONN2}
 Given a manifold\/~$M$ and a linear connection\/ $\nabla$ on\/~$M$ with
 torsion tensor~$T$, we have for any differential form\/~$\alpha$ of
 degree~$r$ and any $r+1$ vector fields $\, X_0^{},\ldots, X_r^{} \,$
 on\/~$M$
 \begin{eqnarray} \label{eq:LEMC1}
  \begin{array}{l}
   {\displaystyle
    \sum_{i=0}^r \, (-1)^i \; (\nabla_{\!X_i}^{} \alpha)
    (X_0^{},\ldots,\hat{X}_i^{},\ldots,X_r^{})} \\
   \quad =~{\displaystyle
            d\alpha(X_0^{},\ldots,X_r^{}) \; +
            \sum_{0 \leqslant i < j \leqslant r}^{\phantom{r}} \, (-1)^{i+j} \;
            \alpha(T(X_i^{},X_j^{}),X_0^{},\ldots,\hat{X}_i^{},\ldots,
                                   \hat{X}_j^{},\ldots,X_r^{})}~.
  \end{array}
 \end{eqnarray}
 Similarly, given a fiber bundle\/~$P$ over a manifold\/~$M$, with vertical
 bundle~$V\!P$, and a partial linear connection\/~$\nabla$ in\/~$V\!P$
 along\/~$V\!P$ with torsion tensor\/~$T$, we have for any vertical
 differential form\/~$\alpha$ of degree~$r$ and any $r+1$ vertical
 vector fields $\, X_0^{},\ldots,X_r^{} \,$ on\/~$P$
 \begin{eqnarray} \label{eq:LEMC2}
  \begin{array}{l}
   {\displaystyle
    \sum_{i=0}^r \, (-1)^i \; (\nabla_{\!X_i}^{} \alpha)
    (X_0^{},\ldots,\hat{X}_i^{},\ldots,X_r^{})} \\
   \quad =~{\displaystyle
            d_V^{} \alpha(X_0^{},\ldots,X_r^{}) \; +
            \sum_{0 \leqslant i < j \leqslant r}^{\phantom{r}} \, (-1)^{i+j} \;
            \alpha(T(X_i^{},X_j^{}),X_0^{},\ldots,\hat{X}_i^{},\ldots,
                                   \hat{X}_j^{},\ldots,X_r^{})}~.
  \end{array}
 \end{eqnarray}
\end{lem}
\proof
 Both statements follow from the same elementary calculation:
 \begin{eqnarray*}
 \lefteqn{\sum_{i=0}^r \, (-1)^i \; (\nabla_{\!X_i}^{} \alpha)
    (X_0^{},\ldots,\hat{X}_i^{},\ldots,X_r^{})} \hspace*{5mm} \\
  &=&\!\! \sum_{i=0}^r \, (-1)^i \; X_i^{} \cdot
          \alpha(X_0^{},\ldots,\hat{X}_i^{},\ldots,X_r^{}) \\
  & & \mbox{} - \, \sum_{i=0}^r \, \sum_{j=0}^{i-1} \, (-1)^i \;
      \alpha \bigl( X_0^{},\ldots,\nabla_{\!X_i}^{} X_j^{},\ldots,
                    \hat{X}_i^{},\ldots,X_r^{} \bigr) \\
  & & \mbox{} - \, \sum_{i=0}^r \, \sum_{j=i+1}^r \, (-1)^i \;
      \alpha \bigl( X_0^{},\ldots,\hat{X}_i^{},\ldots,
                    \nabla_{\!X_i}^{} X_j^{},\ldots,X_r^{} \bigr) \\
  &=&\!\! \sum_{i=0}^r \, (-1)^i \; X_i^{} \cdot
          \alpha(X_0^{},\ldots,\hat{X}_i^{},\ldots,X_r^{}) \\
  & & \mbox{} + \sum_{0 \leqslant i < j \leqslant r}^{\phantom{r}} \, (-1)^{i+j} \;
      \alpha \bigl( \nabla_{\!X_i}^{} X_j^{} - \nabla_{\!X_j}^{} X_i^{} \,,
                    X_0^{},\ldots,\hat{X}_i^{},\ldots,
                    \hat{X}_j^{},\ldots,X_r^{} \bigr) \\[4mm]
  &=&\!\! d\alpha(X_0^{},\ldots,X_r^{}) \quad \mbox{or} \quad
          d_V^{} \alpha(X_0^{},\ldots,X_r^{}) \\[3mm]
  & & \mbox{} + \sum_{0 \leqslant i < j \leqslant r}^{\phantom{r}} \, (-1)^{i+j} \;
      \alpha \bigl( \nabla_{\!X_i}^{} X_j^{} - \nabla_{\!X_j}^{} X_i^{}
                    - [X_i^{},X_j^{}] \,, X_0^{},\ldots,\hat{X}_i^{},\ldots,
                    \hat{X}_j^{},\ldots,X_r^{} \bigr) 
 \end{eqnarray*}
\qed

\section*{Appendix B: Affine manifolds and maps}

We recall some concepts and facts about affine manifolds and affine maps
between them, following~\cite{KN1}.
\begin{dfn}~
 An \textbf{affine manifold} is a manifold equipped with a linear
 connection\/~$\nabla$. A (smooth) map\/ $\, f: M \longrightarrow M' \,$
 between affine manifolds is called an \textbf{affine map} if its
 tangent map\/ $\, Tf: TM \longrightarrow TM' \,$ preserves parallel
 transport, i.e., for any curve\/ $\gamma$ in\/~$M$ from\/ $x$ to\/~$y$
 with image curve\/ $\, \gamma' = f \smcirc \gamma \,$ in\/~$M'$ from\/
 $f(x)$ to\/~$f(y)$, the following diagram commutes:
 \[
  \xymatrix{
   T_x M~ \ar[d]_{U_\gamma^\nabla(x,y)} \ar[r]^{T_x f} &
   ~T_{f(x)} M' \ar[d]^{U_{\gamma'}^{\nabla'}(f(x),f(y))} \\
   T_y M~ \ar[r]_{T_y f} & ~T_{f(y)} M'
  }
 \]
\end{dfn}
Obviously, an affine map takes geodesics into geodesics and therefore
commutes with the exponential, in the sense that
\[
 f \bigl( \exp_x (u) \bigr)~=~\exp_{f(x)} \bigl( T_x f \cdot u \bigr)
 \qquad \mbox{for $\, x \in M$, $u \in \mathrm{Dom}(\exp_x) \subset T_x M$}
\]
(see~\cite[Chapter~6, Proposition~1.1, p.~225]{KN1}).

An important property of riemannian manifolds that extends to affine
manifolds is the existence of convex geodesic balls around each point.
First, we say that an open neighborhood $U_x^{}$ of a point~$x$ in
an affine manifold~$M$ is a \emph{normal neighborhood} of~$x$ if
there is an open neighborhood $U_x^0$ of the origin in~$T_x M$
contained in the domain $\mathrm{Dom}(\exp_x)$ of the exponential
$\exp_x$ such that the latter restricts to a diffeomorphism
$\, \exp_x: U_x^0 \longrightarrow U_x^{}$.
Second, a \emph{geodesic ball} around a point $x$ of~$M$ is a
normal neighborhood $B_x^{}$ of~$x$ obtained as the inverse
image of an open ball in $T_x M$ around the origin, of radius
$\rho$, say, where $\rho$ is sufficiently small, with respect
to some (arbitrarily chosen) scalar product in~$T_x M$.
Now it can be shown~\cite[Chapter~3, Theorem~8.7, p.~149]{KN1}
that geodesic balls $B_x^{}$ of sufficiently small radius~$\rho$
have two additional useful properties: (a)~$B_x^{}$ is geodesically
convex (i.e., any two points of $B_x^{}$ can be connected by a
geodesic entirely contained in~$B_x^{}$) and (b)~$B_x^{}$ is a
normal neighborhood not only of~$x$ but of any of its points.
Whenever this is the case, $B_x^{}$ will be called a
\emph{convex geodesic ball}.
\begin{thm}~ \label{thm:RECOBAF}
 Let\/ $M$ and\/~$M'$ be connected affine manifolds and
 let\/ $\, f: M \longrightarrow M' \,$ be an affine map.
 Suppose that\/ $M$ is simply connected and geodesically 
 complete and that\/ $f$ is a local diffeomorphism.
 Then\/ $f$ is a covering (in particular, it is surjective),
 establishing\/ $M$ as the universal covering manifold of\/~$M'$,
 and\/ $M'$ is also geodesically complete.
\end{thm}
\begin{rmk}~
 The ``riemannian version'' of this theorem (which assumes that $M$
 and $M'$ are riemannian manifolds and $f$ is isometric) is well known
 and can be found in many textbooks, but the proofs given usually make
 use of the Hopf-Rinow theorem and therefore do not extend to the present
 situation, where we do not have metrics (in the topological sense).
 An alternative approach can be found in~\cite[Chapter~10, Theorem~18,
 p.~167]{HIC}, and the proof presented below is an adaptation of that
 to the affine case.
\end{rmk}
\proof
 We begin by showing that $f$ is onto. Considering that $M'$ is connected
 and $f$ is a local diffeomorphism, so that its image $f(M)$ is necessarily
 an open submanifold of~$M'$, it suffices to show that $f(M)$ is also closed.
 Thus let $\, x' \in M' \,$ be a point in the closure of~$f(M)$
 and let $B'$ be a convex geodesic ball in~$M'$ around~$x'$.
 Then there exist a point $\, y' \in B' \cap f(M) \,$ and, due to the
 fact that $B'$ is a normal neighborhood of $y'$ as well, a tangent
 vector $\, u' \in T_{y'} M' \,$ such that $\, \exp_{y'}(u') = x'$.
 Choose $\, y \in M \,$ such that $\, f(y) = y' \,$ and, using
 that $f$ is local diffeomorphism, $u \in T_y M \,$ such that
 $\, T_y f \cdot u = u'$. \linebreak Set $\, x = \exp_y(u)$.
 Then since $f$ is affine, we have $\, f(x) = x'$.
 The argument also shows that $M'$ is geodesically complete.
 Finally, to show that $f$ is a covering, note that the inverse
 image $f^{-1}(x')$ of a point $\, x' \in M' \,$ under the local
 diffeomorphism~$f$ is a discrete subset of~$M$, and we can always
 choose a scalar product on $T_{x'} M'$ and, for every $\, x \in
 f^{-1}(x')$, a scalar product on $T_x M$ such that $\, T_x f:
 T_x M \longrightarrow T_{x'} M' \,$  is isometric; then the
 inverse image under~$f$ of a convex geodesic ball around~$x'$,
 of sufficiently small radius, will be the disjoint union,
 parametrized by $\, x \in f^{-1}(x')$, of the convex geodesic
 balls around~$x$, of the same radius.
\qed

\end{appendix}

\end{document}